\documentclass[10 pt, journal, final, twocolumn]{IEEEtran}
\usepackage{enumerate}
\usepackage{amsmath, amsthm, amssymb}
\usepackage{float}
\usepackage{subcaption}
\usepackage{flexisym}
\usepackage{mathtools}
\usepackage{breqn}
\usepackage{array}
\usepackage{tikz}
\usepackage{cite}
\usetikzlibrary{shapes.geometric, arrows}
\usepackage[font = small]{caption}
\usepackage{paralist}
\usepackage{mathrsfs}
\usepackage{algorithm}
\usepackage{algpseudocode}
\usepackage[inline]{enumitem}
\usepackage{bbm}
\allowdisplaybreaks[1]

\newtheorem{lemma}{Lemma}
\newtheorem{theorem}{Theorem}

\theoremstyle{definition}
\newtheorem{assumption}{Assumption}

\newtheorem{definition}{Definition}

\DeclareMathOperator*{\argmin}{arg\,min}

\DeclareMathOperator*{\diag}{diag}
\DeclareMathOperator*{\blkdiag}{blkdiag}

\renewcommand{\IEEEQED}{\IEEEQEDclosed}

\begin{document}
\title{Resilient Distributed Parameter Estimation with Heterogeneous Data}
\author{Yuan Chen, Soummya Kar, and Jos\'{e} M. F. Moura
\thanks{Yuan Chen {\{(412)-268-7103\}}, Soummya Kar {\{(412)-268-8962\}}, and Jos\'{e} M.F. Moura {\{(412)-268-6341, fax: (412)-268-3890\}} are with the Department of Electrical and Computer Engineering, Carnegie Mellon University, Pittsburgh, PA 15217 {\tt\small \{yuanche1, soummyak, moura\}@andrew.cmu.edu}}
\thanks{This material is based upon work supported by DARPA under agreement number FA8750-16-2-0033, by the Department of Energy under award number DE-OE0000779, and by the National Science Foundation under Award Number CCF 1513936.}}
\maketitle

\begin{abstract}
This paper studies resilient distributed estimation under measurement attacks. A set of agents each makes successive local, linear, noisy measurements of an unknown vector field collected in a vector parameter. The local measurement models are heterogeneous across agents and may be locally unobservable for the unknown parameter. An adversary compromises some of the measurement streams and changes their values arbitrarily. The agents' goal is to cooperate over a peer-to-peer communication network to process their (possibly compromised) local measurements and estimate the value of the unknown vector parameter. We present \textbf{SAGE}, the Saturating Adaptive Gain Estimator, a distributed, recursive, consensus+innovations estimator that is resilient to measurement attacks. We demonstrate that, as long as the number of compromised measurement streams is below a particular bound, then, \textbf{SAGE} guarantees that all of the agents' local estimates converge almost surely to the value of the parameter. The resilience of the estimator -- i.e., the number of compromised measurement streams it can tolerate -- does not depend on the topology of the inter-agent communication network. Finally, we illustrate the performance of \textbf{SAGE} through numerical examples. 
\end{abstract}


\section{Introduction}\label{sect: intro}
The continued growth of the Internet of Things (IoT) has introduced numerous applications featuring networks of devices cooperating toward a common objective. IoT devices instrument smart cities to monitor traffic patterns or air pollution~\cite{AoT}, convoys of driverless cars cooperate to navigate through hazardous weather conditions~\cite{VehiclesExample1, ConnectedVehicles}, and teams of reconnaisance drones form ad-hoc networks to extend their coverage areas~\cite{NetworkedDrones}. These applications require distributed algorithms to locally process measurements collected at each individual device for actionable insight. For example, connected autonomous vehicles need to estimate the global state of traffic from their local sensor measurements for safe and efficient path planning, without transmitting their raw measured data to the cloud. 

IoT devices, however, are vulnerable to malicious cyber attacks, and, in particular, adversaries may pathologically corrupt a subset of the devices' (or agents') measurements~\cite{ChenDistributed2}. Attackers may create false traffic obstacles for a driverless car by spoofing its LIDAR~\cite{VehiclesExample2}, or they may arbitrarily control vehicles by manipulating its onboard sensor measurements~\cite{VehiclesExample3, HACMS}. In IoT setups, attacks on individual devices propagate over the network and affect the behavior of other devices~\cite{ChenIoT, LMSCluster2, ChenDistributed1, ChenDistributed2}. Without proper countermeasures, IoT will be compromised with possibly dangerous consequences.

In this paper, we address the following question: how can an IoT or network of agents resiliently estimate an unknown (vector) parameter while an adversary manipulates a subset of their measurements? Consider a set of agents or devices, connected by a communication network of arbitrary topology, each making a stream of measurements of an unknown parameter $\theta^*$. For example, in air quality monitoring applications, $\theta^*$ represents the field of pollutant concentrations throughout a city, and each device makes successive measurements of a few components of the parameter, corresponding to pollutant concentrations in specific locations. An adversary manipulates a subset of the local measurement streams. The agents' goal is to ensure that \textit{all} agents resiliently estimate the field or (vector) parameter $\theta^*$. This paper presents \textbf{SAGE}, the Saturating Adaptive Gain Estimator, a recursive distributed estimator, that guarantees that, under a sufficient condition on the number of compromised measurement streams, all of the agents' local estimates are (statistically) strongly consistent, i.e., they converge to $\theta^*$ almost surely (a.s.). 

\subsection{Literature Review}
We briefly review related literature, and we start with the Byzantine Generals problem, which demonstrates the effect of adversaries in distributed computation~\cite{ByzantineGenerals}. In this problem, a group of agents must decide whether or not to attack an enemy city by passing messages to each other over an all-to-all communication network. A subset of adversarial agents transmits malicious messages to prevent the remaining agents from reaching the correct decision. Reference~\cite{ByzantineGenerals} shows that a necessary and sufficient condition for \textit{any} distributed protocol to recover the correct decision is that less than one third of the agents be adversarial. Like~\cite{ByzantineGenerals}, reference~\cite{ResilientConsensus} also studies resilient consensus (i.e., reaching agreement among the agents) over all-to-all communication networks. Unlike the consensus problem, which does not incorporate measurements, we study, in this paper, resilient distributed \textit{estimation}, where agents must process their local measurement streams. 

Byzantine attacks also affect performance of decentralized inference algorithms. In decentralized inference, a team of devices measures an unknown parameter and transmits either their measurements or their local decisions to a fusion center, which then processes all data simultaneously to recover the value of the parameter. Adversarial (Byzantine) devices transmit falsified data to the fusion center to disrupt the inference process. References~\cite{Varshney1} and~\cite{Varshney2} address hypothesis testing under Byzantine attack (where the parameter may only take discrete values), and references~\cite{FusionCenterEstimation1} and~\cite{Zhang2} address estimation (where the parameter may take continuous values). The work on resilient decentralized inference in~\cite{Varshney1, Varshney2, FusionCenterEstimation1, Zhang2} restrict the adversary to follow probabilistic attack strategies and use this restriction to design resilient processing algorithms. In this paper, we consider a stronger class of adversaries who may behave \textit{arbitrarily}, unlike the restricted adversaries from~\cite{Varshney1, Varshney2, FusionCenterEstimation1, Zhang2} who must adhere to probabilistic strategies. 

While~\cite{ByzantineGenerals} and~\cite{ResilientConsensus} study resilient computation over fully-connected, all-to-all networks, in general, distributed computation in IoT involves sparse communications, where each device only communicates with a small subset of the other devices. The authors of~\cite{LeBlanc1} propose an algorithm for resilient Byzantine consensus of local agent states over sparse communication networks that guarantees that the non-adversarial agents resiliently reach consensus in the presence of Byzantine adversaries as long as the graph satisfies conditions that restrict topology of the communication network. 
Reference~\cite{LeBlanc2} extends the algorithm in~\cite{LeBlanc1} to address distributed estimation. The estimator in~\cite{LeBlanc2} is resilient to Byzantine agents (i.e., all nonadversarial agents resiliently estimate the unknown parameter) but requires that a known subset of the agents be guaranteed, a priori, to be nonadversarial. In this paper, we do not assume that any agent is impervious to adversaries; any subset of the agents may fall under attack.

A popular method to cope with Byzantine agents is to deploy algorithms that explicitly detect and identify adversaries~\cite{LMSCluster3, Pasqualetti1, Sundaram3, ChenDistributed1}. {\color{black} Reference~\cite{LMSCluster3} studies multi-task distributed estimation, where, instead of estimating a single parameter, the agents may be interested in estimating different parameters. The authors of~\cite{LMSCluster3} propose a distributed algorithm for each agent to identify the neighbors that are interested in estimating the same parameter and ignore information from neighbors that are interested in estimating different parameters. If all adversaries behave as though they measure a different parameter, then, the algorithm from~\cite{LMSCluster3} successfully identifies them. Byzantine agents in general, however, are not required to behave as though they measure and estimate a different parameter; they may send arbitrary information to their neighbors. The algorithm from~\cite{LMSCluster3} is only able to identify one type of Byzantine behavior (i.e., where byzantine agents behave as though they measure a different parameter) but not all Byzantine behaviors in general.}

The authors of~\cite{Pasqualetti1} develop distributed algorithms to detect and identify adversarial agents. In~\cite{Pasqualetti1}, each device must know the structure of the entire network and the consensus algorithm dynamics (i.e., how each agent updates their local state), which is then leveraged to design attack detection and identification filters. {\color{black} In general, attack identification is computationally expensive. Even in centralized settings, where a fusion center has access to all of the agents' measurements, identifying attacks is a combinatorial problem (in the number of measurements)~\cite{Pasqualetti2, Fawzi, Shoukry}. } Reference~\cite{Sundaram3} also requires each agent to know the topology of the entire network, and it designs algorithms to detect and identify adversaries in distributed function calculation. The techniques from~\cite{Pasqualetti1} and~\cite{Sundaram3} become burdensome in terms of memory, communication, and computation requirements as the number of agents and the size of the network grows. {\color{black} To avoid this high computational cost, we develop a resilient distributed estimator that does not depend on expliclitly detecting nor on identifying the attack. We emphasize that not requiring the explicit detection and identification of compromised streams is a feature, not a limitation, of our approach that is very useful in many practical scenarios.}

In our previous work~\cite{ChenDistributed1}, we designed a method to detect the presence of adversaries in distributed estimation that applies to general communication network topologies and only requires each agent to know who its local neighbors are. Our algorithm guaranteed that, in the presence of adversaries, the uncompromised agents will either almost surely recover the value of the unknown parameter or detect the presence of Byzantine agents. Unlike~\cite{ChenDistributed1}, this paper presents a distributed estimation algorithm that copes with attacks on measurement streams \textit{without} explicitly detecting the attack. The estimator in this paper ensures that \textit{all} of the agents, not just the uncompromised agents (like in~\cite{ChenDistributed1}), resiliently recover the value of the unknown parameter. 

Instead of explicitly detecting and identifying attacks as in~\cite{Pasqualetti1, Sundaram3, ChenDistributed1}, other methods cope with adversaries \textit{implicitly}. Implicit countermeasures aim to ensure that the agents achieve their collective objective (e.g., in distributed estimation, the agents' objective is to recover the value of the unknown parameter) without explicitly detecting and identifying the adversary's actions. {\color{black} References~\cite{MultitaskLMS, LMSCluster2} present distributed diffusion algorithms for multi-task estimation, where, like in~\cite{LMSCluster3}, the agents measure and are interested in estimating different parameters. In~\cite{LMSCluster2} and~\cite{MultitaskLMS}, to cope with neighbors who are interested in estimating different parameters, each agent applies adaptive weights to data received from their neighbors. Over time, the agents learn to ignore information that is irrelevant to their estimation tasks. Multi-task distributed estimation algorithms~\cite{MultitaskLMS, LMSCluster2} may cope with adversaries who behave as though they measure a different parameter but fail to cope with Byzantine agents in general. Reference~\cite{ResilientLMS} shows how a single Byzantine agent may prevent distributed diffusion algorithms (\cite{LMSCluster2, LMSCluster3, MultitaskLMS}) from correctly estimating the unknown parameter.}


Our previous work~\cite{ChenDistributed2} presented an algorithm for distributed parameter recovery that we demonstrate analytically to be implicitly resilient against measurement attacks (i.e., attacks that manipulate the agents' local measurement streams). In~\cite{ChenDistributed2}, each agent's measurement was noiseless and homogeneous. That is, in the absence of attacks, every agent had the same, noiseless measurement model. This paper, unlike~\cite{ChenDistributed2}, addresses agents with possibly different or heterogeneous measurement models that are corrupted by noise. A survey and overview of both explicit and implicit security countermeasures for resilient computation in decentralized and distributed IoT architectures is found in our previous work~\cite{ChenIoT}.

\subsection{Summary of Contributions}
This paper studies resilient distributed estimation under measurement attacks. A team of IoT devices (or agents) makes noisy, heterogeneous measurements of an unknown, static parameter $\theta^*$. The devices share information with their neighbors over a communication network and process locally their measurement streams to estimate $\theta^*$. An adversary attacks some of the measurements -- a compromised measurement takes arbitrary value as determined by the adversary. Our goal is to ensure that \textit{all} of the agents, even those with compromised measurements, consistently estimate $\theta^*$. To this end, we present \textbf{SAGE}, the Saturating Adaptive Gain Estimator. 

\textbf{SAGE} is a \textit{consensus+innovations} type estimator where each agent maintains and updates its estimate as a weighted sum of its current estimate, its current measurement, and the estimates of its neighbors  (see, e.g., ~\cite{Kar1, Kar2, Kar3, Kar4}). Each agent applies an adaptive gain to its measurement in the estimate update procedure to limit the impact of aberrant, compromised measurements. We design this gain and establish a sufficient condition on the maximum number of  measurement streams that may be compromised, under which \textbf{SAGE} guarantees that \textit{all} of the agents resiliently estimate the parameter of interest. Specifically, as long as the sufficient condition is satisfied, then, \textbf{SAGE} ensures that all of the agents' estimates are strongly consistent, i.e., the estimate of each (local) agent converges almost surely to the value of the parameter. We will compare this sufficient condition against theoretical bounds of resilience established for \textit{centralized} (or cloud) estimators, where a single processor has access to all of the measurements~\cite{Pasqualetti2, Fawzi, Shoukry}.

{\color{black} Although \textbf{SAGE} shares similar ideas with the Saturated Innovations Update (\textbf{SIU}) algorithm from~\cite{ChenDistributed2}, a resilient distributed parameter recovery algorithm for setups where the agents all make the same (homogeneous), noiseless measurement of the parameter $\theta^*$, \textbf{SAGE} addresses resilient distributed estimation when \begin{enumerate*} \item the agents make different (heterogeneous) measurements and \item the measurements are corrupted by measurement noise. \end{enumerate*} The first extension allows for agents or sensors of different types to cooperate, a more realistic and practical condition. The second makes the analysis of \textbf{SAGE} much more intricate than the analysis of \textbf{SIU}. 

With measurement noise, it is more difficult to distinguish compromised measurements from normal measurements, since a strategic adversary may hide the effect of the measurement attacks within the noise. So, \textbf{SAGE} and the agents must simultaneously cope with the natural disturbance from the noise as well as the artificial disturbance from the attack. With \textbf{SIU} and noiseless settings~\cite{ChenDistributed2}, the agents only have to deal with disturbances resulting from the attack. The analysis of \textbf{SAGE} requires new technical tools, that we develop in this paper, not found in~\cite{ChenDistributed2}, to deal with the measurement noise. Compared to the literature (notably~\cite{ChenDistributed2, Kar1, Kar2, Kar3, Kar4}), the improvements of \textbf{SAGE} are \begin{enumerate*}\item heterogeneous measurement models and measurement noise with respect to~\cite{ChenDistributed2} and \item the presence of attacks with respect to~\cite{Kar1, Kar2, Kar3, Kar4}. \end{enumerate*}} 

This paper focuses on measurement attacks, where an adversary manipulates the values of a subset of the agents' measurement streams, in contrast with Byzantine attacks where adversarial agents send false information to their neighbors. To cope with Byzantine agents, the uncompromised agents must additionally filter and process the information they receive from neighbors~\cite{LeBlanc1, LeBlanc2, Pasqualetti1, LMSCluster2, Sundaram3, ChenDistributed1}. In contrast, to cope with measurement attacks, agents must additionally filter and process their own measurement streams. Existing algorithms for resilient distributed computation with Byzantine agents (e.g.,~\cite{LeBlanc1, LeBlanc2, Pasqualetti1, LMSCluster2, Sundaram3, ChenDistributed1}) only ensure that the uncompromised agents resiliently complete the task. Moreover,  algorithms from~\cite{LeBlanc1, LeBlanc2, Sundaram3} constrain the topology of the communication network to satisfy specific conditions (e.g., in~\cite{Sundaram3}, on the number of unique paths between any two uncompromised agents). In contrast, we will show here that \textbf{SAGE} is resilient to measurement attacks regardless of the network topology, as long as it is connected on average, and we ensure that \textit{all} of the agents, including those with compromised measurements, consistently estimate the unknown parameter. 

The rest of this paper is organized as follows. In Section~\ref{sect: background}, we review the models for the measurements, communications, and attacks. Section~\ref{sect: algorithm} presents the \textbf{SAGE} algorithm, a consensus+innovations estimator that is resilient to measurement attacks. In Section~\ref{sect: converge}, we show that, as long as a sufficient condition on the number of compromised agents is satisfied, then, for any connected network topology, \textbf{SAGE} guarantees that all of the agents' local estimates are strongly consistent. Section~\ref{sect: resilience} compares the sufficient condition (on the number of compromised agents for resilient distributed estimation) against theroetical resilience bounds for centralized (cloud) estimators. We illustrate the performance of \textbf{SAGE} through numerical examples in Section~\ref{sect: examples}, and we conclude in Section~\ref{sect: conclusion}. 
\section{Background}\label{sect: background}
\subsection{Notation}
Let $\mathbb{R}^k$ be the $k$ dimensional Euclidean space, $I_k$ the $k$ by $k$ identity matrix, $\mathbf{1}_k$ and $\mathbf{0}_k$ the $k$ dimensional one and zero vectors. The operator $\left\lVert \cdot \right \rVert_2$ is the $\ell_2$ norm. For matrices $A$ and $B$, $A \otimes B$ is the Kronecker product. A set of vectors $\mathcal{V} = \left\{v_1, \dots, v_{\left\lvert \mathcal{V} \right \rvert} \right\}$ is orthogonal if $v_i ^\intercal v_j = 0$ for all $i \neq j$. 

Let $G = (V, E)$ be a simple undirected graph (no self loops or multiple edges) with vertices $V = \left\{1, \dots, N \right\}$ and edges $E$. The neighborhood of vertex $n$, $\Omega_n$, is the set of vertices that share an edge with $n$. The degree of a node $d_n$ is the size of its neighborhood $\left\lvert \Omega_n \right \rvert$, and the degree matrix of the graph $G$ is $D = \diag (d_1, \dots, d_N)$. The adjacency matrix of $G$ is $A = \left[ A_{nl} \right]$, where $A_{nl} = 1$ if $\left(n, l\right) \in E$ and $A_{nl} = 0$ otherwise. The Laplacian matrix of $G$ is $L = D - A$. The Laplacian matrix $L$ has ordered eigenvalues
$ 0 = \lambda_1(L) \leq \cdots \leq \lambda_N(L),$ and eigenvector $\mathbf{1}_N$ associated with the eigenvalue $\lambda_1(L) = 0.$ For a connected graph $G$, $\lambda_2(L) > 0$. We assume that $G$ is connected in the sequel. References~\cite{Spectral, ModernGraph} review spectral graph theory. 

Random objects are defined on a common probability space $\left( \Omega, \mathcal{F} \right)$ equipped with a filtration $\left\{ \mathcal{F}_t \right\}$. Reference~\cite{Stochastic} reviews stochastic processes and filtrations. Let $\mathbb{P} \left( \cdot \right)$ and $\mathbb{E} \left( \cdot \right)$ be the probability and expectation operators, respectively. In this paper, unless otherwise stated, all inequalities involving random variables hold almost surely (a.s.), i.e., with probability $1$.

\subsection{Measurement and Communication Model}
Consider a set of $N$ agents or devices $\left\{1, 2, \dots, N\right\}$, each making noisy, local streams of measurements of an unknown field represented by a vector parameter $\theta^* \in \mathbb{R}^M$. For example, a network of air quality sensors may monitor several pollutant concentrations over a city. The measurement of agent $n$ at time $t = 0, 1, 2, \dots$ is
\begin{equation}\label{eqn: noAttackMeasurement}
	y_n(t) = H_n \theta^* + w_n(t),
\end{equation}
where $w_n(t)$ is the local measurement noise. Agent $n$ has $P_n$ scalar measurements (at each time $t$), and $M$ is the dimension of the parameter $\theta^*$ (i.e., $H_n \in \mathbb{R}^{P_n \times M}$).{\color{black}\footnote{\color{black}Our measurement model differes from the model used in~\cite{ChenDistributed2}, where each agent made the same, noiseless measurement $y_n(t) = \theta^*$.}}  In air quality monitoring, for example, each device may measure the concentrations of a few pollutants in its neighborhood over time, corresponding to specific components of $\theta^*$.
{\color{black}
\begin{assumption}\label{ass: noise} 
The measurement noise $w_n(t)$ is independently and identically distributed (i.i.d.) over time and independent across agents with mean $\mathbb{E} \left[ w_n(t) \right] = 0$ and covariance $\mathbb{E} \left[w_n(t) w_n(t)^{\intercal} \right] = \Sigma_n.$ The sequence $\left\{ w_n(t) \right\}$ is $\mathcal{F}_{t+1}$ adapted and independent of $\mathcal{F}_t$.
\end{assumption} }

The agents' estimation performance depends on the collection of all of their local measurements. Let
\begin{equation}\label{eqn: stackedMeasurement}
	\mathbf{y}_t = \left[ \begin{array}{ccc} y_1^\intercal (t) & \cdots & y_N^\intercal (t) \end{array} \right]^\intercal = \mathcal{H} \theta^* + \mathbf{w}_t, 
\end{equation}
be stack of all agents' measurements at time $t$, where 
\begin{align*}
	\mathbf{w}_t& = \left[\begin{array}{ccc} w_1^\intercal (t) & \cdots & w_N^\intercal(t) \end{array} \right]^\intercal,\\
	\mathcal{H} &= \left[\begin{array}{ccc} H_1^\intercal & \cdots & H_N^\intercal \end{array} \right]^\intercal
\end{align*}
collect the measurement noises and measurement matrices ($H_1, \dots, H_N$), respectively. Let $P = \sum_{n = 1}^N P_n$ be the total number of (scalar) measurements (at each time $t$) over all agents (i.e., the matrix $\mathcal{H}$ has $P$ rows). {\color{black}We index the scalar components of $\mathbf{y}_t$ from $1$ to $P$, 
\begin{equation} \label{eqn: yRows}
\mathbf{y}_t = \left[ \begin{array}{c} \left. \begin{array}{c} y^{(1)}(t) \\ \vdots \\ y^{(P_1)}(t) \end{array} \right\}y_1(t) \\ \vdots \\ \left. \begin{array}{c} y^{(\overline{P}_n+1)}(t) \\ \vdots \\ y^{(\overline{P}_n + P_n)}(t) \end{array} \right\}y_n(t) \\ \vdots \\ \left. \begin{array}{c} y^{(\overline{P}_N+1)}(t) \\ \vdots \\ y^{(P)}(t) \end{array} \right\}y_N(t) \end{array} \right],
\end{equation} where $\overline{P}_n = \sum_{j=0}^{n-1} P_j$. In general, $\overline{P}_n + 1, \dots, \overline{P}_n + P_n$ index the components of $y_n(t)$. Similar to~\eqref{eqn: yRows},} we label the rows of $\mathcal{H}$ as $h^\intercal_1, \dots, h^\intercal_P$ ($h_p^\intercal \in \mathbb{R}^{1 \times M}$), i.e.,
\begin{equation*}
\mathcal{H} = \left[\begin{array}{ccc} h_1 & \cdots & h_P \end{array} \right]^\intercal,
\end{equation*}
{\color{black} and, in general, $\overline{P}_n + 1, \dots, \overline{P}_n + P_n$, index the rows of $H_n$.}

{\color{black}
\begin{assumption}\label{ass: normalization}
	The measurement vectors $h_p$ are normalized to unit $\ell_2$ norm, i.e., $\left\lVert h_p \right\rVert_2 = 1, \forall p = 1, \dots P$. 
\end{assumption}
\noindent We make Assumption~\ref{ass: normalization} without loss of generality. If Assumption~\ref{ass: normalization} is not satisfied, each agent may compute a normalized version of its local measurement
\begin{equation*}
	\mathcal{D} y_n(t) = \mathcal{D}H_n \theta^* + \mathcal{D} w_n(t),
\end{equation*}
where
\begin{equation*}
	\mathcal{D} = \diag \left(\left\lVert h_{\overline{P_n} + 1} \right\rVert_2^{-1}, \dots, \left\lVert h_{\overline{P_n} + P_n} \right\rVert_2^{-1} \right).
\end{equation*} 
Then, instead of processing the raw measurements $y_n(t)$, the agents may process the normalized measurements $\mathcal{D} y_n(t)$. Note that each row of the matrix $\mathcal{D} H_n$ (i.e., the measurement matrix associated with the normalized measurement  $\mathcal{D} y_n(t)$) satisfies Assumption~\ref{ass: normalization}.}

We now define a measurement stream.
\begin{definition}[Measurement Stream]
A measurement stream $\left\{ y^{(p)}(t) \right\}_{t = 0, 1, 2, \dots}$ is the collection of the scalar measurement $y^{(p)}(t)$ over all times $t = 0, 1, 2, \dots$.\hfill $\small \blacksquare.$
\end{definition}
\noindent For the rest of this paper, we refer to measurement streams by their component index, i.e., measurement stream $p$ refers to $\left\{ y^{(p)}(t) \right\}_{t = 0, 1, 2, \dots}$. Let $\mathcal{P} = \left\{1, \dots, P \right\}$ be the set of all measurement streams.

The agents' collective goal is to estimate the value of $\theta^*$ using their local measurement streams. Global observabilty is an important property for estimation.
\begin{definition}[Global Observability]\label{def: globalObs}
	Let $\mathcal{X} = \left\{p_1, \dots, p_{\left\lvert \mathcal{X} \right\rvert} \right\} \subseteq \mathcal{P}$ be a collection of measurement stream indices. Consider the measurement vectors associated with the streams in $\mathcal{X}$, collected in the stacked measurement matrix
\begin{align*}
	\mathcal{H}_{\mathcal{X}} = \left[\begin{array}{ccc} h_{p_1} & \cdots & h_{p_{\left\lvert \mathcal{X} \right \rvert}} \end{array} \right]^\intercal. \label{eqn: HX}
\end{align*}
The set $\mathcal{X}$ is globally observable if the observability Grammian
\begin{equation}\label{eqn: GX}
	\mathcal{G}_{\mathcal{X}} = \mathcal{H}_{\mathcal{X}}^\intercal \mathcal{H}_{\mathcal{X}} = \sum_{p \in \mathcal{X}} h_p h_p^\intercal
\end{equation}
is invertible. \hfill $\small \blacksquare.$
\end{definition}
\noindent Global observabilty of a set of measurement streams $\mathcal{X}$, means that, in the absence of noise $w_n(t)$, the value of $\theta^*$ may be determined exactly from a single snapshot (in time) of the measurements in $\mathcal{X}$. When the parameter $\theta^*$ and measurement vectors $h_p$ are static, then, having access to streams of measurements over time does not affect observability.

{\color{black}
\begin{assumption}\label{ass: globalObs}
	The set of all measurement streams $\mathcal{P}$ is globally observable: the matrix $\mathcal{G}_{\mathcal{P}} = \sum_{p = 1}^{P} h_p h_p^\intercal$ is invertible.\footnote{Global observability is necessary for centralized estimators to be consistent, so, it is natural to also assume it here for distributed settings. }
\end{assumption}
\noindent In contrast with~\cite{ChenDistributed2}, here, the measurements at each individual agent are not observable, i.e., we do not require local observability, and an individual agent may not consistently estimate $\theta^*$ using just its own local measurements.}

Agents may exchange information with each other over a time-varying communication network, defined by a time varying graph $G(t) = (V, E(t))$. In the graph $G(t)$, the vertex set $V$ is the set of all agents $\left\{1, 2, \dots, N \right\}$, and the edge set $E(t)$ is the set of communication links between agents at time~$t$. 

{\color{black}
\begin{assumption}\label{spec: iidGraphs}
The Laplacians $\left\{L(t)\right\}$ (associated with the graphs $\left\{G(t) \right\}$) form an i.i.d. sequence with mean $\mathbb{E} \left[L(t) \right] = \overline{L}.$ The sequence $\left\{L(t) \right\}$ is $\mathcal{F}_{t+1}$ adapted and independent of $\mathcal{F}_t$. The sequences $\left\{ L(t) \right\}$ and $\left\{w_n(t) \right\}$ are mutually independent. 
\end{assumption}
\noindent The random graphs $G(t)$ (and associated Laplacians $L(t)$) capture communication link failures, like intermittent shadowing common in wireless environments, as well as random gossip communication protocols~\cite{Kar3}. For example, in the absence of link failures, the communication network may be a fixed graph $G = (V, E)$. In practice, however, some of these links may fail intermitently, inducing a random, time varying graph $G(t) = (V, E(t))$. 
}

The agents' collective goal is to estimate the value of the parameter $\theta^*$ (in context of air quality monitoring, the concentration of pollutants) using their measurements and information received from their neighbors.

%

%

\begin{assumption}\label{ass: connectivity}
	The communication network is connected on average. That is, we require $\lambda_2 \left( \overline{L} \right) > 0$. We do \textit{not} require $G(t)$ to be connected. 
\end{assumption} 

\subsection{Attack Model}
An attacker attempts to disrupt the agents and prevent them from estimating $\theta^*$. The attacker replaces a subset of the measurements with arbitrary values. We model the effect of the attack as an additive disturbance $a_n(t)$, so that, under attack, the measurement of agent $n$ is
\begin{equation}\label{eqn: attackMeasurement}
	y_n(t) = H_n \theta^* + w_n(t) + a_n(t).
\end{equation}
Similar to $y_n(t)$ and $H_n$, we label the components of $a_n(t)$ by $\overline{P}_n + 1, \dots, \overline{P}_n + P_n$, i.e.,
	$a_n(t) = \left[\begin{array}{ccc} a^{\left( \overline{P}_n + 1 \right)} (t) & \cdots & a^{\left( \overline{P}_n + P_n \right)} (t)\end{array} \right]^\intercal.$
We say that measurement stream $p$ is compromised or under attack if, \textit{for any} $t = 0, 1, \dots,$ $a^{(p)}(t) \neq 0$.\footnote{A compromised measurement stream does not have to be under attack for \textit{all} times $t$. If there is any $t$ such that $a^{(p)}(t) \neq 0$, we consider $p$ to be compromised even if there exists $t' \neq t$ such that $a^{(p)}(t) = 0.$} {\color{black} We make no assumptions about \textit{how} the attacker compromises the measurement streams. The attacker may choose $a_n(t)$ arbitrarily, which, in general, may be unbounded and time varying. }For example, in air quality monitoring, if measurement stream $p$ is compromised, then, at some time $t$, instead of taking the value of the true pollutant concentration, $y^{(p)}(t)$ becomes an arbitrarily valued concentration measurement. 
For analysis, we partition $\mathcal{P}$, the set of all measurement streams, into the set of compromised measurement streams \begin{equation} \mathcal{A} = \left\{p \in \mathcal{P} \Big\vert \exists t = 0, 1, \dots, \: a^{(p)}(t) \neq 0 \right\}, \end{equation} and the set of noncompromised or uncompromised measurement streams $\mathcal{N} = \mathcal{P} \setminus \mathcal{A}$. Each agent $n$ may have both compromised and uncompromised local measurement streams: if one of agent $n$'s measurement streams is compromised, it does not necessarily mean that all of its streams are compromised. The sets $\mathcal{A}$ and $\mathcal{N}$ are time invariant. Agents do not know which measurement streams are compromised. 


When measurements fall under attack, resilient estimation of $\theta^*$ depends on the notion of (global) $s$-sparse observabilty.

\begin{definition}[Global $s$-Sparse Observabilty (\!\cite{Fawzi, Shoukry})] Let $\mathcal{P}$ be the set of \textit{all} measurement streams. The set of all measurement streams is globally $s$-sparse observable if the matrix
\begin{equation*}
	\mathcal{G}_{\mathcal{P} \setminus \mathcal{X}} = \mathcal{H}^\intercal_{\mathcal{P}\setminus \mathcal{X}} \mathcal{H}_{\mathcal{P}\setminus \mathcal{X}} = \sum_{p \in \mathcal{P}\setminus \mathcal{X}} h_p h_p^\intercal.
\end{equation*}
is invertible for \textit{all} subsets $\mathcal{X} \subseteq \mathcal{P}$ of cardinality $\left\lvert \mathcal{X} \right \rvert = s$. \hfill $\small \blacksquare.$
\end{definition}
\noindent A set of measurements is (globally) $s$-sparse observable if it is still (globally) observable after removing \textit{any} $s$ measurements from the set~\cite{Fawzi, Shoukry}. 

We make the following assumptions on the attacker:
\begin{assumption}\label{ass: attackSetSize}
	The attacker may only manipulate a subset of the measurement streams. That is, $0 \leq \frac{\left\lvert \mathcal{A} \right \rvert}{P} < 1$. The attacker may not manipulate \textit{all} of the measurement streams.
\end{assumption}
\begin{assumption}\label{ass: measurementAttacksOnly}
	The attack does not change the value of $\theta^*$. That is the attacker may only manipulate the agents' measurements of $\theta^*$, but it may not change the true value of $\theta^*$. In the context of air quality monitoring, this means that the attacker may change some of the  \textit{measurements} of the pollutant concentrations, but may not change the true concentrations themselves.
\end{assumption}

The attack model in this paper differs from the Byzantine attack model (see, e.g.,~\cite{ByzantineGenerals, ChenDistributed1, Varshney1, Sundaram1, LeBlanc1, LeBlanc2}). Under the Byzantine attack model, a subset of adversarial agents behaves arbitrarily and transmits false information to their neighbors, while the remaining nonadversarial agents must cope with pathologically corrupted data received from their adversarial neighbors. In contrast, under the measurement attack model~\eqref{eqn: attackMeasurement}, agents do not intentionally transmit false information to their neighbors. Instead, agents must cope with adversarial data from a subset of their measurement streams. Resilient algorithms for the Byzantine attacker model aim to ensure that only the non-adversarial, uncompromised agents accomplish their computation task~\cite{ChenDistributed1, LeBlanc1, LeBlanc2}. This is because Byzantine agents may behave arbitrarily and deviate from any prescribed protocol. For distributed estimation under the measurement attack model (equation~\eqref{eqn: attackMeasurement} with Assumptions~\ref{ass: attackSetSize} and~\ref{ass: measurementAttacksOnly}), we aim to ensure that \textit{all} agents recover the value of the unknown parameter, even those agents with compromised measurement streams. 

{\color{black} Byzantine attacks compromise the communication between agents. Attackers may intercept communication packets and change their contents, or they may, via malicious software, force agents to transmit false information to neighbors~\cite{ZhangIoT}. The measurement attack model in \textbf{SAGE} considers attacks that compromise the data collected by the agents. Attackers spoof the agents' measurements so that they collect false data. For example, attackers may create false obstacles for a driverless car by feeding it false LIDAR measurements~\cite{VehiclesExample2}. }

\section{\textbf{SAGE}: Saturating Adaptive Gain Estimator}\label{sect: algorithm}
In this section, we present \textbf{SAGE}, the Saturating Adaptive Gain Estimator, a distributed estimation algorithm that is resilient to measurement attacks. $\textbf{SAGE}$ is a \textit{consensus + innovations} distributed estimator~\cite{Kar1, Kar2, Kar3, Kar4} and is based on the Saturated Innovations Update (\textbf{SIU}) algorithm we introduced~\cite{ChenDistributed2}. Unlike the \textbf{SIU} algorithm, which requires the agents' measurement streams to be noiseless (i.e., $w_n(t) \equiv 0$ in~\eqref{eqn: noAttackMeasurement} and~\eqref{eqn: attackMeasurement}) and homogeneous ($H_n \equiv I_M$), \textbf{SAGE} may simultaneously cope with disturbances from both measurement noise and the attacker on heterogeneous (possibly different $H_n$) linear measurement models. 

\subsection{Algorithm Description}
\textbf{SAGE} is an iterative algorithm in which each agent maintains and updates a local estimate $x_n(t)$ based on its local observation $y_n(t)$ and the information it receives from its neighbors. Each iteration of \textbf{SAGE} consists of a message passing step and an estimate update step. Each agent $n$ sets its initial local estimate as $x_n(0) = 0$. 

\underline{Message Passing}: Every agent $n$ transmits its current estimate, $x_n(t)$, to each of its neighbors $l \in \Omega_n$. 

\underline{Estimate Update}: Every agent $n$ maintains a running average of its local measurement stream; i.e., each agent computes
\begin{equation}\label{eqn: avgMeasurement}
	\overline{y}_n(t) = \frac{t}{t+1} \overline{y}_n(t-1) + \frac{1}{t+1} y_n(t),
\end{equation}
with initial condition $\overline{y}_n(-1) = 0$. Note that, if $p$ is not under attack (i.e., $p \in \mathcal{N}$), then, 
\begin{equation*}
	\overline{y}^{(p)}(t) = h_{ p}^\intercal \theta^* + \overline{w}^{(p)}(t),
\end{equation*}
where $\overline{w}^{(p)}_n(t) = \frac{1}{t+1} \sum_{j = 0}^t w^{(p)}_n(j)$ is time-averaged measurement noise. Every agent $n$ updates its estimate as
\begin{equation}\label{eqn: estimateUpdate}
\begin{split}
	x_n(t+1) =& x_n(t) - \beta_t \sum_{l \in \Omega_n(t)} \left(x_n(t) - x_l(t) \right) \\
		& + \alpha_t H_n^\intercal K_n(t)  \left(\overline{y}_n(t) - H_n x_n(t) \right),
\end{split}
\end{equation}
where $\alpha_t > 0$ and $\beta_t > 0$ are weight sequences to be defined in the sequel. The term $\overline{y}_n(t) - H_n x_n(t)$ is \textit{innovation} of agent $n$. The innovation term is the difference between the agent's observed (averaged) measurement $\overline{y}_n(t)$ and its predicted measurement $H_n x_n(t)$, based on its current estimate $x_n(t)$ and the measurement model~\eqref{eqn: noAttackMeasurement}. The innovation term incorporates information from the agent's measurement streams into its updated estimate. The second term in~\eqref{eqn: estimateUpdate} is the global consensus average. 

The term $K_n(t)$ is a diagonal matrix of gains that depends on the innovation $\overline{y}_n(t) - H_n x_n(t)$ and is defined as
\begin{equation}\label{eqn: KnDef}
	K_n(t) = \diag \left( k_{\overline{P}_n+ 1}(t), \dots k_{\overline{P}_n + P_n}(t) \right),
\end{equation}
where, for $p = \overline{P}_n+1, \dots, \overline{P}_n + P_n$,
\begin{equation}\label{eqn: KnDef2}
	k_{p}(t) = \min \left(1, \gamma_t{\left\lvert \overline{y}^{\left(p\right)}_n(t) - h_{p}^\intercal x_n(t) \right\rvert}^{-1} \right),
\end{equation}
and $\gamma_t$ is a scalar adaptive threshold to be defined in the sequel. The gain $k_p(t)$ saturates the innovation term, along each component, at the threshold $\gamma_t$. That is, $k_p(t)$ ensures that, for each measurement stream $p$, the magnitude of the scaled innovation, $\left \lvert k_p(t) \left(\overline{y}_n^{(p)}(t) - h_p^\intercal x_n(t) \right) \right\rvert$, never exceeds the threshold $\gamma_t$. 
%
%
%

The purpose of $K_n(t)$ and $\gamma_t$ is to limit the impact of attacks. It may, however, also limit the impact of measurements that are not under attack ($p \in \mathcal{N}$). The main challenge in designing the \textbf{SAGE} algorithm is to choose the adaptive threshold $\gamma_t$ to balance these two effects. We choose $\gamma_t$ and the weight sequences $\alpha_t$ and $\beta_t$ as follows:
\begin{enumerate}
\item Select $\alpha_t$, $\beta_t$ as:
\begin{equation}\label{eqn: alphaBeta}
	\alpha_t = \frac{a}{(t+1)^{\tau_1}}, \beta_t = \frac{b}{(t+1)^{\tau_2}},
\end{equation}
with $a, b > 0$ and $0 < \tau_2 < \tau_1 < 1$.
\item Select $\gamma_t$ as:
\begin{equation}\label{eqn: gammaDef}
	\gamma_t = \frac{\Gamma} {(t+1)^{\tau_{\gamma}}},
\end{equation}
with $\Gamma > 0$ and $0 < \tau_\gamma < \min\left(\frac{1}{2} , \tau_1 - \tau_2 \right).$ 
\end{enumerate}

{\color{black} Choosing $\alpha_t$, $\beta_t$, and $\gamma_t$ correctly is key to the resilience and consistency of \textbf{SAGE}. If we choose these weight sequences incorrectly, then, \textbf{SAGE} is not resilient to the measurement attacks. We carefully craft $\alpha_t$ and $\beta_t$ to decay over time, with $\alpha_t$ decaying faster than $\beta_t$. In the long run, the consensus term dominates the innovation term. The difference in decay rates results in the mixed time-scale behavior of the estimate update~\eqref{eqn: estimateUpdate}, which is crucial for consistency~\cite{Kar1, Kar2, Kar3, Kar4}.}
The threshold $\gamma_t$ also decays over time. Initially, the estimate update~\eqref{eqn: estimateUpdate}, with large contributions from the innovation terms, incorporates more information from the measurements. After several updates, as local estimates $x_n(t)$ become closer to the true value of the parameter, we expect the innovation to decrease in magnitude for uncompromised measurements. By reducing the threshold $\gamma_t$ over time, we limit the impact of compromised measurements without overly affecting the contribution from the noncompromised measurements.

{\color{black} Compared with \textbf{SIU}~\cite{ChenDistributed2}, \textbf{SAGE} uses the same innovation and consensus weights $\alpha_t$ and $\beta_t$~\eqref{eqn: alphaBeta} but uses a different threshold sequence $\gamma_t$~\eqref{eqn: gammaDef}. In~\cite{ChenDistributed2}, which considered homogeneous measurements without noise, we designed threshold $\gamma_t$ to be an upper bound on the magnitude of the innovation term for noncompromised measurements. That is, in~\cite{ChenDistributed2}, if measurement stream $p$ is not under attack, then, the innovation component associated with $p$ is guaranteed to have magnitude less than $\gamma_t$. When there is (unbounded) measurement noise, it is impossible to construct a threshold sequence with this property: even if measurement stream $p$ is not under attack, the noise may cause the innovation magnitude to exceed the threshold. So, instead of using the threshold sequence from~\cite{ChenDistributed2}, we use a new threshold sequence~\eqref{eqn: gammaDef}, which is more suitable for the noisy setting. Further, in \textbf{SAGE}, because the local measurements are heterogeneous, they can also, and will in general be, unobservable (that we refer to in this paper as locally unobservable). This differs with \textbf{SIU} where the agents measure directly the unknown (vector) parameter and so, in the absence of attacks, are (locally) observable. }


\subsection{Main Results: Strong Consistency of \textbf{SAGE}}
We present our main results on the performance of \textbf{SAGE}. The first result establishes the almost sure convergence of the \textbf{SAGE} algorithm (under a condition on the measurement streams under attack) and specifies its rate of convergence. Recall that $\mathcal{A}$ is the set of compromised measurement streams.
\begin{theorem}[\textbf{SAGE} Strong Consistency]\label{thm: main}
	Let $\mathcal{A} = \left\{p_1, \dots, p_{\left\lvert \mathcal{A} \right \rvert} \right\}$, $\mathcal{H}_{\mathcal{A}} = \left[ \begin{array} {ccc} h_{p_1} & \cdots & h_{p_{\left\lvert \mathcal{A} \right\rvert}} \end{array} \right]^\intercal$, and $\mathcal{N} = \mathcal{P} \setminus \mathcal{A}$. If the matrix 
\begin{equation}\label{eqn: normalGrammian}
	\mathcal{G}_{\mathcal{N}} = \sum_{p \in \mathcal{N}} h_p h_p^\intercal,
\end{equation}
satisfies
\begin{equation}\label{eqn: resilienceCondition}
	\lambda_\min \left(\mathcal{G}_{\mathcal{N}} \right) > \Delta_{\mathcal{A}},
\end{equation}
where
\begin{equation}\label{eqn: maxDisturbance}
	 \Delta_{\mathcal{A}} =  \max \limits_{v \in \mathbb{R}^{\left\lvert \mathcal{A} \right\rvert}, \left\lVert v \right \rVert_\infty \leq 1} \left\lVert \mathcal{H}_{\mathcal{A}}^\intercal v \right\rVert_2,
\end{equation} then, for all $n$, we have
\begin{equation}\label{eqn: localConsistency}
	\mathbb{P} \left( \lim_{t \rightarrow \infty} \left(t+1 \right)^{\tau_0} \left\lVert x_n(t) - \theta^* \right\rVert_2 = 0 \right) = 1,
\end{equation}
for every $0 \leq \tau_0 < \min\left(\tau_\gamma, \frac{1}{2} - \tau_\gamma\right).$
\end{theorem}
{\color{black} We are interested in determining the attacks against which \textbf{SAGE} is resilient. Theorem~\ref{thm: main} provides a sufficient condition for the attacks on the measurements under which \textbf{SAGE} ensures strong consistency of all local estimates.\footnote{\color{black}The ceiling at $\frac{1}{2}$ on the rate of convergence $\tau_0$ in  Theorem~\ref{thm: main} comes as a result of the measurement noise. This is the fastest rate at which the norm of the averaged measurement noise converges to $0$. This ceiling is not present in~\cite{ChenDistributed2} because it assumed noiseless measurements.} Due to heterogeneous measurements, this resilience condition~\eqref{eqn: resilienceCondition} is different from the resilience condition in~\cite{ChenDistributed2} (which assumed homogeneous measurements). This is necessary because, heterogeneous measurements provide different information about $\theta^*$ (e.g., when each measurement stream provides information about a single component of $\theta^*$), so the resilience of \textbf{SAGE} depends on \textit{which} measurement streams are compromised. 

This condition depends on $\mathcal{G}_{\mathcal{N}}$~\eqref{eqn: normalGrammian}, the observability Grammian of the noncompromised measurements, and $\Delta_{\mathcal{A}}$~\eqref{eqn: maxDisturbance}, which captures the maximum disturbance caused by the compromised measurements. The minmum eigenvalue of $\mathcal{G}_{\mathcal{N}}$, $\lambda_{\min} \left( \mathcal{G}_{\mathcal{N}} \right)$ captures the redundancy of the noncompromised measurement streams in measuring $\theta^*$. For example, consider the case where all measurement vectors $h_p$ are canonical basis vectors, so that each scalar measurement stream provides information about one component of $\theta^*$. Then, there are at least $\lambda_{\min} \left( \mathcal{G}_{\mathcal{N}} \right)$ scalar streams measuring each component of~$\theta^*$. 

Intuitively, the resilience condition~\eqref{eqn: resilienceCondition} states that there needs to be enough redundancy in the noncompromised measurements $\mathcal{N}$ to overcome the effects of any attack on the measurement streams in $\mathcal{A}$. Condition~\eqref{eqn: resilienceCondition} is sufficient, but it is not always necessary for consistent estimation. As we will demonstrate through numerical examples, in some cases (i.e., for certain sets of compromised measurements $\mathcal{A}$), \textbf{SAGE} still guarantees strongly consistent local estimates even if~\eqref{eqn: resilienceCondition} is not satisfied.}



The optimization problem in~\eqref{eqn: maxDisturbance} is nonconvex (since it is a constrained maximization of a convex function). As a consequence of Assumption~\ref{ass: normalization}, which states that $\left\lVert h_p \right\rVert_2 = 1$, we have $\left\lvert \mathcal{A} \right \rvert \geq \Delta_{\mathcal{A}}$. Then, a sufficient condition for~\eqref{eqn: resilienceCondition} is
\begin{equation}\label{eqn: relaxedResilience}
	\lambda_\min \left( \mathcal{G}_{\mathcal{N}} \right) > \left\lvert \mathcal{A} \right \rvert.
\end{equation}
The condition in~\eqref{eqn: relaxedResilience} describes the total \textit{number} of compromised measurement streams that are tolerable under \textbf{SAGE}, and we are interested in comparing this number against theoretical limits established in~\cite{Shoukry} for centralized (cloud) estimators. 

Specifically, Theorem 3.2 in~\cite{Shoukry} states that, in the absence of measurement noise, a centralized estimator may tolerate attacks on \textit{any} $s$ measurement streams and still consistently recover the value of the parameter if and only if the set of all measurement streams is $2s$-sparse observable. {\color{black} By tolerating attacks on any $s$ measurement streams, we mean that, no matter \textit{which} $s$ measurement streams are under attack, and no matter \textit{how} the attacker maniuplates those measurements, the estimator still produces consistent estimates. If an estimator is not resilient to attacks on $s$ measurements, it means that there is a specific set of $s$ measurements, and a specific way in which they are manipulated, such that the resulting estimate is not consistent. If the attacker compromises a different set of $s$ measurements, or changes the way the measurements are manipulated, the estimator may still produce consistent estimates. We again emphasize that the resilience condition~\eqref{eqn: relaxedResilience} is a sufficient condition; under certain attacks, even if~\eqref{eqn: relaxedResilience} is not satisfied, \textbf{SAGE} still produces consistently local estimates.} 

The following result relates the (relaxed) resilience condition~\eqref{eqn: relaxedResilience} to sparse observability. 
\begin{theorem}\label{thm: resilience}
	If $\lambda_\min\left(\mathcal{G}_{\mathcal{N}} \right) > \left\lvert \mathcal{A} \right \rvert$ for all $\left\lvert \mathcal{A} \right \rvert \leq s$, then, the set of all measurement streams $\mathcal{P} = \left\{1, \dots, P \right\}$ is globally $2s$-sparse observable. Now, let
\begin{equation}\label{eqn: uniqueH}
	\mathcal{V} = \bigcup_{p \in \mathcal {P}} \left\{ h_p \right\},
\end{equation}
be the set of all unique measurement vectors $h_p$. If $\mathcal{V}$ is orthogonal, then $\lambda_\min\left(\mathcal{G}_{\mathcal{N}} \right) > \left\lvert \mathcal{A} \right \rvert$ only if $\mathcal{P}$ is $2s$-sparse observable.

\end{theorem}

\noindent In general, $2s$-sparse observability is a necessary condition for \textbf{SAGE} to tolerate attacks on \textit{any} $s$ sensors. If the set of all \textit{unique} measurement vectors $h_p$ is orthogonal, then, $2s$-sparse observability is also sufficient. 

The scenario where $\mathcal{V}$, the set of all unique measurement vectors, is orthogonal arises when we are interested in estimating a high dimensional field parameter. For example, in air pollution monitoring, each component of the parameter $\theta^*$, which may be high dimensional, represents the pollutant concentration at a particular location. Each individual device is only able to measure concentrations of nearby locations. In this case, the rows $h_p^\intercal$ of the measurement matrix $H_n$ are the canonical basis (row) vectors, where each device measures a subset of the components of $\theta^*$. {\color{black} Recall that agents are interested in estimating the entire parameter $\theta^*$. When the measurement vectors $h_p$ are canonical basis vectors, the local measurements at each agent provide information about some, but not all, components of $\theta^*$. In order to estimate all of the components of $\theta^*$ when the measurements are not locally observable, the agents must exchange information with their neighbors.}

{\color{black} For analyzing the resilience condition $\lambda_{\min} \left(\mathcal{G}_{\mathcal{N}} \right) > \left\lvert \mathcal{A} \right \rvert$, another case of interest is when $\theta^*$ is a scalar parameter. When $\theta^*$ is a scalar parameter, each (uncompromised) agent's measurement stream is modeled by
\begin{equation}\label{eqn: scalarMeasurement}
	y_n(t) = \theta^* + a_n(t). 
\end{equation}
For the measurement model~\eqref{eqn: scalarMeasurement}, the resilience condition $\lambda_{\min}\left(\mathcal{G}_{\mathcal{N}} \right) > \left\lvert \mathcal{A} \right \rvert$ becomes
$\frac{\left\lvert \mathcal{A} \right \rvert}{P} < \frac{1}{2}. $
To ensure that \textbf{SAGE} produces strongly consistent local estimates of scalar parameters, it is sufficient to have a majority of the measurement streams remain \textit{uncompromised}. For the scalar parameter and measurement model~\eqref{eqn: scalarMeasurement}, the condition $\frac{\left\lvert \mathcal{A} \right \rvert}{P} < \frac{1}{2}$ is also necessary. If a majority of the measurement streams are compromised, then, the adversary may force the local estimates to converge to any arbitrary value. This is because the adversary may compromise measurement streams arbitrarily; $a_n(t)$ may take any (i.e., unbounded) value in~\eqref{eqn: attackMeasurement}. The condition $\frac{\left\lvert \mathcal{A} \right \rvert}{P} < \frac{1}{2}$ is also necessary for resilient \textit{centralized} estimation~\cite{Shoukry}.}

\section{Convergence Analysis: Proof of Theorem~\ref{thm: main}}\label{sect: converge}
This section proves Theorem~\ref{thm: main}, which states that, under \textbf{SAGE}, all local estimates are strongly consistent as long as the resilience condition~\eqref{eqn: resilienceCondition} is satisfied ($\lambda_{\min} \left(\mathcal{G}_{\mathcal{N}} \right) > \Delta_{\mathcal{A}}$).  The proof of Theorem~\ref{thm: main} requires several intermediate results that are found in the Appendix. The proof of Theorem~\ref{thm: main} has two main components: we separately show that, in the presence of measurement attacks, the local estimates $x_n(t)$ converge a.s. to the mean of all of the agents' estimates (i.e., the agents reach consensus on their local estimates), and that, under resilience condition~\eqref{eqn: resilienceCondition}, the mean estimate converges a.s. to the true value of the parameter. Unless otherwise stated, all inequalities involving random variables hold a.s. (with probability $1$). 

\subsection{Resilient Consensus}
In this subsection, we show that the agents reach consensus. 
Let
\begin{align}
	\mathbf{x}_t &= \left[\begin{array}{ccc} x_1(t)^\intercal & \cdots & x_N(t)^\intercal \end{array} \right]^\intercal, \\
	\overline{\mathbf{y}}_t &= \left[\begin{array}{ccc} \overline{y}_1(t)^\intercal & \cdots & \overline{y}_N(t)^\intercal \end{array} \right]^\intercal, \\
	\mathbf{K}_t &= \diag \left(k_1(t), \dots, k_P(t) \right),\\
	D_H &= \blkdiag \left(H_1, \dots, H_N\right), \\
	\overline{\mathbf{x}}_t &= \frac{1}{N} \left( \mathbf{1}^\intercal_N \otimes I_M \right) \mathbf{x}_t, \\
	\widehat{\mathbf{x}}_t & = \mathbf{x}_t - \left(\mathbf{1}_N \otimes I_M \right)\overline{\mathbf{x}}_t.	
\end{align}
The vectors $\mathbf{x}_t, \overline{\mathbf{y}}_t \in \mathbb{R}^{NM}$ collect the estimates and time-averaged measurements, respectively, of all of the agents. The matrix $\mathbf{K}_t \in \mathbb{R}^{P \times P}$ is a diagonal matrix whose diagonal elements are the saturating gains $k_{p}(t)$ for each measurement stream. The variable $\overline{\mathbf{x}}_t \in \mathbb{R}^M$ is the network average estimate (the mean of the agents' local estimates) and $\widehat{\mathbf{x}} \in \mathbb{R}^{NM}$ collects the differences between the agents' local estimate $x_n(t)$ and the network average estimate $\overline{\mathbf{x}}_t$.

\begin{lemma}\label{lem: resilientConsensus}
	Under \textbf{SAGE}, $\widehat{\mathbf{x}}$ satisfies
\begin{equation}\label{eqn: resilientConsensus}
	\mathbb{P} \left(\lim_{t \rightarrow \infty} (t+1)^{\tau_3} \left\lVert \widehat{\mathbf{x}}_t \right \rVert_2 = 0 \right) = 1,
\end{equation}
for every $0 \leq \tau_3 < \tau_\gamma + \tau_1 - \tau_2.$
\end{lemma}

Lemma~\ref{lem: resilientConsensus} states that, under \textbf{SAGE}, the agents' local estimates converge a.s. to the network average estimate at a polynomial rate. That is, the agents reach consensus on their local estimates even in the presence of measurement attacks. Under \textbf{SAGE}, the agents reaching consensus do \textit{not} depend on the number of compromised measurement streams: the result in Lemma~\ref{lem: resilientConsensus} holds even when \textit{all} of the measurement streams fall under attack, as long as the inter-agent communication network is connected on average. 

\begin{IEEEproof}[Proof of Lemma~\ref{lem: resilientConsensus}]
	From~\eqref{eqn: estimateUpdate}, we find that $\mathbf{x}_t$ follows the dynamics
\begin{equation}\label{eqn: stackedDynamics}
	\begin{split}
		\mathbf{x}_{t+1} &= \left( I_{NM} - \beta_t L_t \otimes I_M - \alpha_t D_H^\intercal \mathbf{K}_t D_H \right) \mathbf{x}_t + \\
		& \quad \alpha_t D_H^\intercal \mathbf{K}_t \overline{\mathbf{y}}_t.
	\end{split}
\end{equation}
Then, from~\eqref{eqn: stackedDynamics}, we have that $\widehat{\mathbf{x}}_t$ follows the dynamics
\begin{equation}\label{eqn: resCon1}
\begin{split}
	\widehat{\mathbf{x}}_{t+1} &= \left( I_{NM} -P_{NM} - \beta_t L_t \otimes I_M \right) \widehat{\mathbf{x}}_t + \\
	& \quad \alpha_t \left(I_{NM} - P_{NM} \right) D_H^\intercal \mathbf{K}_t  \left( \overline{\mathbf{y}}_t - D_H\mathbf{x}_t \right),
\end{split}
\end{equation}
where
\begin{equation}\label{eqn: PNMDef}
	P_{NM} = \frac{1}{N}\left(\mathbf{1}_N\mathbf{1}_N\right)^\intercal \otimes I_M.
\end{equation}
From~\cite{Kar1}, we have that the eigenvalues of the matrix $I_{NM} - P_{NM} - \beta_t L \otimes I_M$ are $0$ and $1 - \lambda_n \left(L \right)$, for $n = 2, \dots, N$, each repeated $M$ times. 

Taking the $\ell_2$ norm of $\widehat{\mathbf{x}}_{t+1}$, we have
\begin{equation}\label{eqn: resCon2}
\begin{split}
	\!\!\left\lVert \widehat{\mathbf{x}}_{t+1} \right \rVert_2\!\! &\leq\!\! \left\lVert \left( I_{NM}\! -\! \beta_t L_t \!\otimes\!I_M \!-\! P_{NM} \right) \widehat{\mathbf{x}}_t \right\rVert_2 + \alpha_t \gamma_t \sqrt{P},
\end{split}
\end{equation}
where we have used the following relations to derive~\eqref{eqn: resCon2} from~\eqref{eqn: resCon1}:
\begin{align}
	\left\lVert I_{NM} - P_{NM} \right \rVert_2  &= 1, \label{eqn: resCon3} \\
	\left\lVert h_p \right \rVert_2 &= 1, \label{eqn: resCon3a} \\
	\left\lVert K_n(t) \left(\overline{y}_n(t) - H_n x_n(t) \right) \right \rVert_\infty &\leq \gamma_t. \label{eqn: resCon4}
\end{align}
Relation~\eqref{eqn: resCon3} follows from the definitions of $I_{NM}$ and $P_{NM}$, relation~\eqref{eqn: resCon3a} follows from Assumption~\ref{ass: normalization}, and relation~\eqref{eqn: resCon4} follows from the definition of $K_n(t)$ (equation~\eqref{eqn: KnDef}). By definition, $\widehat{\mathbf{x}}_t \in \mathcal{C}^{\perp}$, where $\mathcal{C} \subseteq \mathbb{R}^{NM}$ is the consensus subspace
\begin{equation}
	\mathcal{C} = \left\{ w \in \mathbb{R}^{NM} \vert w = \mathbf{1}_N \otimes v, v \in \mathbb{R}^N \right\},
\end{equation}
and $\mathcal{C}^\perp$ is the orthogonal complement of $\mathcal{C}$. Since $\widehat{\mathbf{x}}_t \in \mathcal{C}^\perp$, we have $P_{NM} \widehat{\mathbf{x}} = 0$. Lemma~\ref{lem: randomGraphWeights} in the appendix provides an upper bound for the first term on the right hand side of~\eqref{eqn: resCon2}. Specifically Lemma~\ref{lem: randomGraphWeights} states that, for $t$ large enough, we have
\begin{equation}\label{eqn: resCon5}
	\left\lVert \left(I_{NM} - \beta_t L_t \otimes I_M\right) \widehat{\mathbf{x}}_t \right\rVert_2 \leq \left(1 - r(t) \right) \left\lVert \widehat{\mathbf{x}}_t \right\rVert_2,
\end{equation}
where $r(t)$ is an $\mathcal{F}_{t+1}$ adapted process that satisfies $0 \leq r(t) \leq 1$ a.s. and $\mathbb{E} \left[ r(t) \vert \mathcal{F}_t \right] \geq \frac{c_r}{(t+1)^{\tau_2}}$ for some $c_r > 0$. Substituting~\eqref{eqn: resCon5} into~\eqref{eqn: resCon2}, we have
\begin{equation}\label{eqn: resCon6}
	\left\lVert \widehat{\mathbf{x}}_{t+1} \right \rVert_2 \leq \left(1 - r(t) \right) \left\lVert \widehat{\mathbf{x}}_t \right \rVert_2 + \alpha_t \gamma_t \sqrt{P}.
\end{equation}
Lemma~\ref{lem: timeVaryingStochasticConvergence} in the appendix studies the stability of the relation in~\eqref{eqn: resCon6} and shows that $(t+1)^{\tau_3}\left\lVert \widehat{\mathbf{x}}\right\rVert_2$ converges a.s. to $0$ for every $0 \leq \tau_3 < \tau_\gamma + \tau_1 - \tau_2$.
\end{IEEEproof}

\subsection{Network Average Behavior}
This subsection characterizes the behavior of the network average estimate $\overline{\mathbf{x}}_t$. Specifically, we show that, under resilience condition~\eqref{eqn: resilienceCondition}, $\overline{\mathbf{x}}_t$ converges a.s. to the value of the parameter $\theta^*$. Define
\begin{align}
	\overline{\mathbf{e}}_t &= \overline{\mathbf{x}}_t - \theta^*,\label{eqn: networkAverageError} \\
	\widetilde{k}_p(t)&= \left\{\begin{array}{ll} k_p(t), & \text{if } p \in \mathcal{N}, \\
		0, & \text{otherwise,} \end{array} \right. \\
	\mathbf{K}^{\mathcal{N}}_t &= \diag \left(\widetilde{k}_1(t), \dots, \widetilde{k}_p(t) \right), \\
	\mathbf{K}^{\mathcal{A}}_t &= \mathbf{K}_t - \mathbf{K}^{\mathcal{N}}_t, \\
	\overline{\mathbf{w}}_t &= \left[ \begin{array}{ccc} \overline{w}_1(t)^\intercal & \cdots & \overline{w}_N(t)^\intercal \end{array} \right]^\intercal. 
\end{align}
The term $\overline{\mathbf{e}}_t$ in~\eqref{eqn: networkAverageError} is the network average estimation error.
\begin{lemma}\label{lem: staysBounded}
Define the auxiliary threshold
\begin{equation}\label{eqn: lineGammaDef}
	\overline{\gamma}_t = \frac{\Gamma}{(t+1)^{\tau_\gamma}} - \frac{Y}{(t+1)^{\tau_3}} - \frac{W}{(t+1)^{\frac{1}{2} - \epsilon_W}},
\end{equation}
where \begin{equation}\label{eqn: tau3Def}\tau_3 = \tau_\gamma + \tau_1 - \tau_2 - \epsilon_Y\end{equation} for arbitrarily small \begin{equation}\label{eqn: epsilonDef}0 < \epsilon_Y < \tau_1 - \tau_2, \: 0 < \epsilon_W < \frac{1}{2} - \tau_\gamma.\end{equation}
If the resilence condition~\eqref{eqn: resilienceCondition} is satisfied, then, almost surely, there exists $T_0 \geq 0, 0< Y < \infty, 0 < W < \infty$ such that, 
\begin{enumerate}
\item $\left\lVert \overline{\mathbf{w}}_t \right\rVert_2 \leq \frac{W}{(t+1)^{\frac{1}{2} - \epsilon_W}}$ a.s.,
\item $\left\lVert \widehat{\mathbf{x}}_t \right\rVert_2 \leq \frac{Y}{(t+1)^{\tau_3} }$ a.s., and,
\item if, for some $T_1 \geq T_0$, we have $\left\lVert \overline{\mathbf{e}}_{T_1} \right\rVert_2 \leq \overline{\gamma}_{T_1}$, then, $\left\lVert \overline{\mathbf{e}}_t \right \rVert_2 \leq \overline{\gamma}_t$ for all $t \geq T_1$ a.s.
\end{enumerate}
\end{lemma}

\begin{IEEEproof}[Proof of Lemma~\ref{lem: staysBounded}]
	From~\eqref{eqn: estimateUpdate}, we find that $\mathbf{e}_t$ follows the dynamics
\begin{equation}\label{eqn: averageErrorDynamics}
	\begin{split}
		\overline{\mathbf{e}}_{t+1} &= \left( I_{M} - \frac{\alpha_t}{N} \sum_{p \in \mathcal{N}} k_p(t) h_p h_p^\intercal \right) \overline{\mathbf{e}}_t - \\
		& \quad \frac{\alpha_t}{N} \left( \mathbf{1}_N^\intercal \otimes I_M \right) D_H^\intercal \mathbf{K}^{\mathcal{N}}_t \left(D_H\widehat{\mathbf{x}}_t - \overline{\mathbf{w}}_t \right)+ \\
	&\quad \frac{\alpha_t}{N} \left( \mathbf{1}_N^\intercal \otimes I_M \right) D_H^\intercal \mathbf{K}^{\mathcal{A}}_t \left( {\overline{\mathbf{y}}}_t - D_H \mathbf{x}_t \right).
	\end{split}
\end{equation}
Lemma~\ref{lem: measureNoise} in the appendix characterizes the a.s. convergence of the time-averaged noise $\overline{\mathbf{w}}_t$ and shows that
\begin{equation}\label{eqn: staysBounded1}
	\mathbb{P} \left(\lim_{t \rightarrow \infty} (t+1)^{\delta_0} \left\lVert \overline{\mathbf{w}}_t \right\rVert_2 = 0 \right) = 1,
\end{equation}
for every $0 \leq \delta_0 < \frac{1}{2}.$ As a result of Lemma~\ref{lem: resilientConsensus}, we have
\begin{equation}\label{eqn: staysBounded2}
	\mathbb{P} \left(\lim_{t \rightarrow \infty} (t+1) ^{\tau_3} \left \lVert \widehat{\mathbf{x}}_t \right \rVert_2 = 0 \right) = 1,
\end{equation}
for every $0 \leq \tau_3 < \tau_\gamma + \tau_1 - \tau_2$. 

We now consider the evolution $\overline{\mathbf{e}}_t$ along sample paths $\omega \in \Omega$ such that 
\begin{align}
	\lim_{t \rightarrow \infty} (t+1)^{\delta_0} \left\lVert \overline{\mathbf{w}}_{t, \omega}\right\rVert_2 &= 0, \label{eqn: staysBounded3} \\
	\lim_{t \rightarrow \infty} (t+1) ^{\tau_3} \left \lVert \widehat{\mathbf{x}}_{t, \omega} \right \rVert_2 &= 0. \label{eqn: staysBounded4}
\end{align}
Note that the set of all such $\omega$ has probability measure $1$. The notation $\overline{\mathbf{e}}_{t, \omega}$ means the value of $\overline{\mathbf{e}}_t$ along the sample path~$\omega$. As a result of~\eqref{eqn: staysBounded3} and~\eqref{eqn: staysBounded4}, there exists finite $Y_\omega$ and $W_\omega$ such that, with probability $1$, 
\begin{align}
	\left\lVert \overline{\mathbf{w}}_{t, \omega} \right \rVert_2 &\leq \frac{W_\omega} {(t+1)^{\frac{1}{2} - \epsilon_W}}, \: \left\lVert \widehat{\mathbf{x}}_{t, \omega} \right \rVert_2 \leq \frac{Y_\omega} {(t+1)^{\tau_3}}. 
\end{align}
Suppose that, at time $T_1$, $\left\lVert \overline{\mathbf{e}}_{T_1, \omega} \right\rVert_2 \leq \overline{\gamma}_{T_1, \omega}$ (where, in~\eqref{eqn: lineGammaDef}, $Y = Y_\omega$ and $W = W_\omega$). We show that there exists finite $T_0$ such that if $T_1 \geq T_0$, then, for all $t \geq T_1$, $\left\lVert \mathbf{e}_{t, \omega} \right \rVert \leq \overline{\gamma}_{t, \omega}$ a.s.

Specifically, we show that, for sufficiently large $T_1$, if $\left\lVert \overline{\mathbf{e}}_{T_1, \omega} \right \rVert_2 \leq \overline{\gamma}_{T_1, \omega}$ then, $\left\lVert \overline{\mathbf{e}}_{T_1 + 1, \omega} \right \rVert_2 \leq \overline{\gamma}_{T_1 + 1, \omega}$. {\color{black} We show that  if $\left\lVert \overline{\mathbf{e}}_{T_1, \omega} \right \rVert_2 \leq \overline{\gamma}_{T_1, \omega}$, then, $k_p(T_1, \omega) = 1$ for all $p \in \mathcal{N}$ (i.e., the magnitude of the innovations for noncompromised measurements $p$ are all below the threshold $\gamma_t$.) } By the triangle inequality, we have, for all $p \in \mathcal{N}$,
\begin{align}
\begin{split}
	&\left\lvert\overline{y}^{(p)}_n(T_1, \omega) - h_p^\intercal x_n(T_1, \omega) \right \rvert \\
	&\: \leq \left\lvert h_p^\intercal \left( \theta^* - \overline{\mathbf{x}}_{T_1, \omega} + \overline{\mathbf{x}}_{T_1, \omega} - x_n(T_1, \omega) \right) \right\rvert + \\ &\quad \:\:  \left\lvert \overline{w}^{p}_n(T_1, \omega) \right \rvert, \\
\end{split} \\
\begin{split}
	&\: \leq \left\lVert \overline{\mathbf{e}}_{T_1, \omega} \right \rVert_2 + \left\lVert \widehat{\mathbf{x}}_{T_1, \omega} \right\rVert_2 + \left\lVert \overline{\mathbf{w}}_{T_1, \omega} \right \rVert_2,
\end{split}\\
&\: \leq \overline{\gamma}_{T_1, \omega} + \frac{Y_\omega}{(T_1+1)^{\tau_3}} + \frac{W_\omega}{(T_1+1)^{\frac{1}{2} - \epsilon_W}}.
\end{align}
By definition of $\overline{\gamma}_t$, we have
\begin{equation}\label{eqn: staysBounded8}
	\overline{\gamma}_{t, \omega} = \gamma_t -  \frac{Y_\omega}{(t+1)^{\tau_3}} - \frac{W_\omega}{(t+1)^{\frac{1}{2} - \epsilon_W}},
\end{equation}
which means that
\begin{equation}\label{eqn: staysBounded9}
	\left\lvert y^{(p)}_n(T_1, \omega) - h_p^\intercal x_n(T_1, \omega) \right \rvert \leq \gamma_{T_1}, \forall p \in \mathcal{N}. 
\end{equation}
As a result of~\eqref{eqn: staysBounded9}, we have $k_p(T_1, \omega) = 1$ for all $p \in \mathcal{N}$.

Using~\eqref{eqn: averageErrorDynamics}, we analyze the dynamics of $\left\lvert \overline{\mathbf{e}}_t \right\rVert_2$. From~\eqref{eqn: averageErrorDynamics}, we have
\begin{equation}\label{eqn: staysBounded7}
	\begin{split}
		\left\lVert \overline{\mathbf{e}}_{t+1, \omega} \right\rVert_2 & \leq \left \lVert I_{M} - \frac{\alpha_t}{N} \sum_{p \in \mathcal{N}} k_p(t) h_p h_p^\intercal \right \rVert_2 \left \lVert \overline{\mathbf{e}}_{t, \omega} \right \rVert_2 + \\
 	& \quad \alpha_t \frac{\left\lvert \mathcal{N} \right \rvert}{N} \left( \left\lVert \widehat{\mathbf{x}}_{t, \omega} \right \rVert_2 + \left \lVert \overline{\mathbf{w}}_{t, \omega} \right \rVert_2\right) + \alpha_t \frac{\Delta_{\mathcal{A}}}{N} \gamma_t.
	\end{split}
\end{equation}
To derive~\eqref{eqn: staysBounded7}, we use the fact that $k_p(t) \leq 1$ and $ k_p(t) \left\lvert \overline{y}^{(p)}_n(t) - h_p^\intercal (t) x_n(t) \right\rvert \leq \gamma_t$ for all $p \in \mathcal{A}$ and for all~$t$. {\color{black}Substituting~\eqref{eqn: staysBounded8} into~\eqref{eqn: staysBounded7}, and using the fact that  $k_p(T_1, \omega) = 1$ for all $p \in \mathcal{N}$,} we have
\begin{equation}\label{eqn: staysBounded10}
\begin{split}
	\left\lVert \overline{\mathbf{e}}_{T_1 + 1, \omega} \right \rVert_2 &\leq \left(1 -\frac{\alpha_{T_1}\kappa}{N} \right) \overline{\gamma}_{T_1, \omega} + \frac{\alpha_{T_1} P Y_\omega}{N(T_1+1)^{\tau_3}} +\\ &\quad  \frac{\alpha_{T_1} P W_\omega}{N (T_1+1)^{\frac{1}{2} - \epsilon_W}},
\end{split}
\end{equation} 
where \begin{equation} \kappa = \lambda_\min \left( \mathcal{G} \right) - \Delta_{\mathcal{A}}. \end{equation} To derive~\eqref{eqn: staysBounded10}, we also use the fact that, for $T_1$ sufficiently large, $\frac{\alpha_{T_1} \lambda_\min \left( \mathcal{G}_{\mathcal{N}} \right)}{N} \leq 1$ and
\begin{equation}\label{eqn: normCondition}
	\left\lVert I_M - \frac{\alpha_{T_1}}{N}\sum_{p \in \mathcal N} h_p h_p^\intercal \right\rVert_2 = 1 - \frac{\alpha_{T_1}}{N} \lambda_\min \left( \mathcal{G}_{\mathcal{N}} \right).
\end{equation}
Equation~\eqref{eqn: normCondition} follows from the fact that, by definition, $\mathcal{G}_{\mathcal{N}} = \sum_{p \in \mathcal{N}} h_p h_p^\intercal$ (equation~\eqref{eqn: normalGrammian}). 

{\color{black} Using~\eqref{eqn: staysBounded10}, we now find conditions on $T_1$ such that $\left\lVert \overline{\mathbf{e}}_{T_1 + 1, \omega} \right \rVert_2 \leq \overline{\gamma}_{T_1 + 1, \omega}$.} Define
\begin{equation}
	\overline{\Gamma}_{t, \omega} = \Gamma - \frac{Y_\omega}{(t+1)^{\tau_3 - \tau_\gamma}} - \frac{W_\omega}{ (t+1)^{\frac{1}{2} - \epsilon_W - \tau_{\gamma} }}.
\end{equation}
Note that $\overline{\gamma}_{t, \omega} = \frac{\overline{\Gamma}_{t, \omega}} {(t+1)^{\tau_\gamma}}$, $\overline{\Gamma}_{t, \omega}$ is increasing in $t$, and, for $t$ sufficiently large, $\overline{\Gamma}_{t, \omega} > 0$. Without loss of generality, let $\widetilde{T}_\omega$ be defined such that, for all $t \geq \widetilde{T}_\omega$,~\eqref{eqn: normCondition} holds and $\overline{\Gamma}_{t, \omega} > 0$. Since $\overline{\Gamma}_{t, \omega}$ is increasing in $t$, we have
\begin{equation}\label{eqn: staysBounded11}
	\overline{\gamma}_{t+1, \omega} \geq \frac{\overline{\Gamma}_{t, \omega}}{(t+2)^{\tau_\gamma}} = \left(\frac{t+1}{t+2} \right)^{\tau_\gamma} \overline{\gamma}_{t, \omega}. 
\end{equation}
Thus, to find conditions on $T_1$ such that \[\left\lVert \overline{\mathbf{e}}_{T_1 + 1, \omega} \right \rVert_2 \leq \overline{\gamma}_{T_1 + 1, \omega},\] it is sufficient to find conditions on $T_1$ such that 
\begin{equation}\label{eqn: staysBounded11}
	\left\lVert \overline{\mathbf{e}}_{T_1 + 1, \omega} \right \rVert_2 \leq \left(\frac{T_1+1}{T_1+2} \right)^{\tau_\gamma} \overline{\gamma}_{T_1, \omega}. 
\end{equation}

Since~$\left\lVert \overline{\mathbf{e}}_{T_1, \omega} \right \rVert_2 \leq \overline{\gamma}_{T_1, \omega}$, we may rearrange the right hand side of~\eqref{eqn: staysBounded10} so that~\eqref{eqn: staysBounded11} becomes
\begin{equation}
\begin{split}\label{eqn: staysBounded12}
	\left(1 - \frac{\alpha_{T_1} \rho_{T_1, \omega}}{N} \right) \overline{\gamma}_{T_1, \omega} \leq \left(\frac{T_1+1}{T_1+2} \right)^{\tau_\gamma} \overline{\gamma}_{T_1, \omega},
\end{split}
\end{equation}
where
\begin{equation}\label{eqn: staysBounded13}
\begin{split}
	\rho_{T_1, \omega} &= \kappa - \frac{P Y_{\omega}}{\overline{\Gamma}_{T_1, \omega}(T_1+1)^{\tau_3 - \tau_\gamma}} - \\ &\quad \frac{P W_{\omega}}{\overline{\Gamma}_{T_1, \omega} (T_1+1)^{\frac{1}{2} - \epsilon_W - \tau_\gamma}}.
\end{split}
\end{equation}
Note that
\begin{equation}\label{eqn: staysBounded14}
	\frac{\alpha_{T_1} \rho_{T_1, \omega}}{N} \geq \frac{\alpha_{T_1}}{N} \left(\kappa - \frac{c_5}{(T_1 + 1)^{\tau_4 - \tau_\gamma}} \right),
\end{equation} 
where $\tau_4 = \min \left(\tau_3, \frac{1}{2} - \epsilon_W \right)$ and $c_5 = \frac{P \left(Y_\omega + W_\omega\right)}{\overline{\Gamma}_{T_1, \omega}}$. For \[T_1 \geq \max \left( \widetilde{T}_\omega, \left(\frac{c_5}{\kappa}\right)^{\frac{1}{\tau_4 - \tau_\gamma}} -1 \right),\] the right hand side of~\eqref{eqn: staysBounded14} is nonnegative, and, after substitution of~\eqref{eqn: staysBounded14}, the condition in~\eqref{eqn: staysBounded12} becomes
\begin{equation}\label{eqn: staysBounded16}
\begin{split}
	1 - &\left(\frac{a\kappa N^{-1}}{(T_1 + 1)^{\tau_1}} - \frac{a c_5 N^{-1}}{(T_1 + 1)^{\tau_1 + \tau_4 - \tau_\gamma}} \right) \leq \left( \frac{T_1 + 1}{T_1 + 2} \right)^{\!\!\tau_\gamma}\!\!,\!\!
\end{split}
\end{equation}
where we have used the fact that, by definition (equation~\eqref{eqn: alphaBeta}), $\alpha_{T_1} = \frac{a}{(T_1+1)^{\tau_1}}$. 

Since $(1-x) \leq e^{-x}$ for $x \geq 0$,~\eqref{eqn: staysBounded16} becomes
\begin{equation}\label{eqn: staysBounded17}
	\begin{split}
		& \frac{a c_5 (N\tau_\gamma)^{-1}}{(T_1 +1)^{\tau_1 + \tau_4 - \tau_\gamma}} -\frac{a\kappa (N\tau_\gamma)^{-1}}{(T_1 + 1)^{\tau_1}}\leq \log \left(\frac{T_1 + 1} {T_1 + 2} \right).
	\end{split}
\end{equation}
Note that $\log \left(\frac{T_1 + 1} {T_1 + 2} \right) \geq 1 -\frac{T_1 + 2} {T_1 + 1} = -\frac{1}{T_1 + 1}$, so a sufficient condition for~\eqref{eqn: staysBounded17} is
\begin{equation}\label{eqn: staysBounded18}
	 \frac{a c_5}{(T_1 +1)^{\tau_1 + \tau_4 - \tau_\gamma}}-\frac{a\kappa}{(T_1 + 1)^{\tau_1}} \leq \frac{-N\tau_\gamma}{T_1 + 1}. 
\end{equation}
Let $\tau_5 = \min(\tau_1 + \tau_4 - \tau_\gamma, 1)$. By definition, $\tau_5 > \tau_1$. Then,~\eqref{eqn: staysBounded18} is satisfied for all $T_1 \geq \left(\frac{a c_5 + N\tau_\gamma}{a \kappa}\right)^{\frac{1}{\tau_5 - \tau_1}} - 1$. Define
\begin{equation}\label{eqn: staysBounded19}
\begin{split}
	&T_0\!\! =\!\! \max\left( \widetilde{T}_\omega, \left(\frac{c_5}{\kappa}\right)^{\frac{1}{\tau_4 - \tau_\gamma}}, \left(\frac{a c_5 + N \tau_\gamma}{a \kappa}\right)^{\frac{1}{\tau_5 - \tau_1}} \right) \!\!-\!\! 1.
\end{split}
\end{equation}
If $T_1 \geq T_0$ and $\left\lVert \overline{\mathbf{e}}_{T_1, \omega} \right \rVert_2 \leq \overline{\gamma}_{T_1, \omega}$, then $\left\lVert \overline{\mathbf{e}}_{T_1 + 1, \omega} \right \rVert_2 \leq \overline{\gamma}_{T_1 + 1, \omega}$. Since $T_1 + 1 \geq T_0$, we then have $\left\lVert \overline{\mathbf{e}}_{t, \omega} \right \rVert_2 \leq \overline{\gamma}_{t, \omega}$ for all $t \geq T_1$. This analysis holds for almost all sample paths $\omega$, and Lemma~\ref{lem: staysBounded} follows. 
\end{IEEEproof}

Lemma~\ref{lem: staysBounded} states that, with probability $1$, there exists some finite time $T_0$ such that, if the $\ell_2$-norm of the network average error at some time $T_1 \geq T_0$ falls below $\overline{\gamma}_{T_1}$, then, it will be upper bounded by $\overline{\gamma}_{t}$ for all $t \geq T_1$. The auxiliary threshold $\overline{\gamma}_t$ is defined in such a way that, if $\left\lVert \overline{\mathbf{e}}_t \right \rVert_2 \leq \overline{\gamma}_t$, then, with probability $1$, $k_p(t) = 1$ for all $p \in \mathcal{N}$. That is, if the network average error has $\ell_2$ norm less than the auxiliary threshold $\gamma_t$, then, \textbf{SAGE} does not affect the contribution of the uncompromised measurements in the agents' local estimate updates. We now use Lemma~\ref{lem: staysBounded} above to show that the network average estimate converges almost surely to $\theta^*$. 

\begin{lemma}\label{lem: averageConv}
	Under \textbf{SAGE}, if $\lambda_{\min}\left(\mathcal{G}_{\mathcal{N}} \right) > \Delta_{\mathcal{A}},$ then
\begin{equation}\label{eqn: averageConv}
	\mathbb{P} \left( \lim_{t \rightarrow \infty} (t+1)^{\tau_0} \left\lVert \overline{\mathbf{e}}_t \right \rVert_2 = 0 \right) = 1,
\end{equation}
for every $0 \leq \tau_0 < \min\left(\tau_\gamma, \frac{1}{2} - \tau_\gamma \right)$. 
\end{lemma}

\begin{IEEEproof}[Proof of Lemma~\ref{lem: averageConv}]
	We study the convergence of $\overline{\mathbf{e}}_t$ along sample paths $\omega \in \Omega$ for which Lemma~\ref{lem: staysBounded} holds (such a set of sample paths has probabilty measure $1$). Consider a specific instantiation of $\left\lVert \overline{\mathbf{e}}_{t, \omega} \right \rVert_2$. There exists finite $T_0$ such that, if at any time $T_1 \geq T_0$, $\left\lVert \overline{\mathbf{e}}_{T_1, \omega} \right \rVert_2 \leq \overline{\gamma}_{T_1, \omega}$, then, for all $t \leq T_1$,  $\left\lVert \overline{\mathbf{e}}_{t, \omega} \right \rVert_2 \leq \overline{\gamma}_{t, \omega}$. Consider $t > T_0$. There are two possibilities: either 1.) there exists $T_1 > T_0$ such that $\left\lVert \overline{\mathbf{e}}_{T_1, \omega} \right \rVert_2 \leq \overline{\gamma}_{T_1, \omega}$, or 2.) for all $t \geq T_0$, $\left\lVert \overline{\mathbf{e}}_{t, \omega} \right \rVert_2 > \overline{\gamma}_{t, \omega}.$ If the first case occurs, then, by Lemma~\ref{lem: staysBounded}, for $t \geq T_1$,
\begin{equation}\label{eqn: avgConv1}
	\left\lVert \overline{\mathbf{e}}_{t, \omega} \right \rVert_2 \leq \frac{\Gamma}{(t+1)^{\tau_\gamma}},
\end{equation}
which means that $\lim_{t \rightarrow \infty} (t+1)^{\tau_0} \left\lVert \mathbf{e}_t \right \rVert_2 = 0$ for every $0 \leq \tau_0 < \tau_\gamma$. 

If the second case occurs, then, define
\begin{equation}\label{eqn: avgConv2}
	\widehat{K}_{t, \omega} = \frac{\overline{\gamma}_{t, \omega} + \frac{Y_\omega}{(t+1)^{\tau_3}} + \frac{W_\omega}{(t+1)^{\frac{1}{2} - \epsilon_W}}}{\left\lVert \overline{\mathbf{e}}_{t, \omega} \right \rVert_2  + \frac{Y_\omega}{(t+1)^{\tau_3}} + \frac{W_\omega}{(t+1)^{\frac{1}{2} - \epsilon_W}}}.
\end{equation}
Since $\left\lVert \overline{\mathbf{e}}_{t, \omega} \right \rVert_2 > \overline{\gamma}_{t, \omega}$, $\widehat{K}_{t, \omega} < 1$. By definition, the numerator in~\eqref{eqn: avgConv2} is equal to $\gamma_t$. Recall that the sample $\omega$ is chosen from a set (of probabilty measure $1$) on which Lemma~\ref{lem: staysBounded} holds, which means that there exists finite $Y_\omega$ and $W_\omega$ that satisfy $\left\lVert \widehat{\mathbf{x}}_{t, \omega} \right\rVert_2 \leq \frac{Y_\omega}{(t+1)^{\tau_3}}$ and $\left\lVert \overline{\mathbf{w}}_{t, \omega} \right \rVert_2 \leq \frac{W_\omega}{(t+1)^{\frac{1}{2} - \epsilon_W}}$. As a result, the denominator in~\eqref{eqn: avgConv2} satisfies
\[\left\lVert \overline{\mathbf{e}}_{t, \omega} \right \rVert_2  + \frac{Y_\omega}{(t+1)^{\tau_3}} + \frac{W_\omega}{(t+1)^{\frac{1}{2} - \epsilon_W}} \geq \left\lvert \overline{y}^{(p)}_n(t, \omega) - x^{(p)}_n(t, \omega) \right\rvert, \]
for all $p \in \mathcal{N}$. Thus, we have
\begin{equation}\label{eqn: avgConv3}
	\widehat{K}_{t, \omega} < k_p(t, \omega),
\end{equation}
for all $n \in \mathcal{N}$. Rearranging~\eqref{eqn: avgConv2}, we also have
\begin{equation}\label{eqn: avgConv4}
	\gamma_t = \widehat{K}_{t, \omega} \left(\left\lVert \overline{\mathbf{e}}_{t, \omega} \right \rVert_2  + \frac{Y_\omega}{(t+1)^{\tau_3}} + \frac{W_\omega}{(t+1)^{\frac{1}{2} - \epsilon_W}} \right).
\end{equation}

Substituting~\eqref{eqn: avgConv3} and~\eqref{eqn: avgConv4} into~\eqref{eqn: staysBounded7}, we have
\begin{align}
\begin{split}\label{eqn: avgConv5}
	&\left\lVert \overline{\mathbf{e}}_{t + 1, \omega} \right \rVert_2 \leq \left(1 -\frac{\alpha_{t}}{N}\widehat{K}_{t, \omega}\lambda_\min\left(\mathcal{G}_{\mathcal{N}} \right)\right) \left\lVert \overline{\mathbf{e}}_{t, \omega} \right \rVert_2 + \\ &\quad \alpha_{t} \frac{\left\lvert \mathcal{N} \right \rvert}{N} \left(\frac{Y_\omega}{(t+1)^{\tau_3}} + \frac{W_\omega}{(t+1)^{\frac{1}{2} - \epsilon_W}}\right) + \alpha_t \frac{\Delta_{\mathcal{A}}}{N} \gamma_t,
\end{split}\\
\begin{split} \label{eqn: avgConv6}
	&= \left(1 - \frac{\alpha_t}{N} \kappa \widehat{K}_{t, \omega} \right) \left\lVert \overline{\mathbf{e}}_{t, \omega} \right \rVert_2 + \\
	& \quad \frac{\alpha_{t}P}{N} \left(\frac{Y_\omega}{(t+1)^{\tau_3}} + \frac{W_\omega}{(t+1)^{\frac{1}{2} - \epsilon_W}}\right).
\end{split}
\end{align}
To derive~\eqref{eqn: avgConv5} from~\eqref{eqn: staysBounded7}, we have used the fact that, for $t$ large enough, 
\begin{align}
	\left \lVert 1 - \frac{\alpha_t}{N} \sum_{p \in \mathcal{N}} k_p(t, \omega) h_p h_p^\intercal \right\rVert_2 &\leq \left \lVert 1 - \frac{\alpha_t}{N} \widehat{K}_{t, \omega} \mathcal{G}_{\mathcal{N}} \right\rVert_2, 
\end{align}
since, by definition, $\mathcal{G}_{\mathcal{N}} = \sum_{p \in \mathcal{N}} h_p h_p^\intercal.$
From~\eqref{eqn: avgConv6}, we show that $\left(t+1\right)^{\tau_\gamma}\widehat{K}_{t, \omega} \geq \widehat{K}_\omega$ for some $\widehat{K}_\omega > 0$. It suffices to show that $\left\lVert \overline{\mathbf{e}}_{t, \omega} \right \rVert < \infty$.  Note that there exists finite $T$ such that $\alpha_T \leq 1$. Because such a finite $T$ exists, to show $\left\lVert \overline{\mathbf{e}}_{t, \omega} \right \rVert < \infty$, it then suffices to only consider $t \geq T$ (since, for $t < T$, $\left\lVert \overline{\mathbf{e}}_{t, \omega} \right \rVert_2 < \infty$ by definition). From~\eqref{eqn: avgConv2}, we have
\begin{equation}\label{eqn: avgConv7}
	(t+1)^{\tau_\gamma} \widehat{K}_{t, \omega} \geq \frac{\Gamma}{\sup_{j \in \left[T, t\right]} \left\lVert \overline{\mathbf{e}}_{j, \omega} \right \rVert_2 + R_{\omega}},
\end{equation}
where $R_\omega = Y_\omega + W_\omega.$

Define the system
\begin{equation}\label{eqn: avgConv8}
	\begin{split}
	\overline{m}_{t+1} &= \left(1 - \frac{\alpha_t}{N} \frac{\kappa \gamma_t}{m_t + R_\omega} \right) m_t + \frac{\alpha_t P}{N} \frac{R_\omega}{(t+1)^{\tau_6}}, \\
	m_{t+1} &= \max\left(\overline{m}_{t+1}, m_t \right),
	\end{split}
\end{equation}
where $\tau_6 = \min\left( \tau_3, \frac{1}{2} -\epsilon_W \right)$, with initial condition $m_T = \left\lVert \overline{\mathbf{e}}_{T, \omega} \right \rVert_2$. By definition, $m_t \geq \sup_{j \in \left[t, T\right]} \left\lVert \overline{\mathbf{e}}_{j, \omega} \right \rVert_2$. Further, from~\eqref{eqn: epsilonDef}, we have $\tau_\gamma < \tau_6$, since $\tau_\gamma < \frac{1}{2} - \epsilon_W$ and $\tau_\gamma < \tau_3 = \tau_\gamma + \tau_1 - \tau_2$\footnote{By definition~\eqref{eqn: alphaBeta}, $\tau_1 > \tau_2$.}, so the system in~\eqref{eqn: avgConv8} falls under the purview of Lemma~\ref{lem: supBounded}. As a result, we have
	$\sup_{t \geq T} \left\lVert \overline{\mathbf{e}}_{t, \omega} \right\rVert_2 < \infty.$
Since $ \left\lVert \overline{\mathbf{e}}_{t, \omega} \right\rVert_2$ is bounded from above, we have that
\begin{equation}\label{eqn: avgConv10}
	\widehat{K}_{t, \omega} \geq \frac{\widehat{K}_\omega}{(t+1)^{\tau_\gamma}},
\end{equation}
for some $\widehat{K}_\omega > 0$. Substituting~\eqref{eqn: avgConv10} into~\eqref{eqn: avgConv6}, we have
\begin{equation}\label{eqn: avgConv11}
\begin{split} 
	 \left\lVert \overline{\mathbf{e}}_{t+1, \omega} \right \rVert_2&\leq \left(1 -  \frac{\alpha_t\kappa \widehat{K}_\omega}{N(t+1)^{\tau_\gamma}} \right) \left\lVert \overline{\mathbf{e}}_{t, \omega} \right \rVert_2 + \frac{\alpha_{t}R_\omega}{(t+1)^{\tau_6}}.
\end{split}
\end{equation}
Lemma~\ref{lem: timeVaryingConvergence} in the appendix studies the convergence of the systems in the form of~\eqref{eqn: avgConv11}. As a result, we have
\begin{equation}\label{eqn: avgConv12}
	\lim_{t \rightarrow \infty} \left(t+1\right)^{\tau_0} \left\lVert \overline{\mathbf{e}}_{t, \omega} \right\rVert_2 = 0,
\end{equation}
for every $0 \leq \tau_0 < \tau_6 - \tau_\gamma$. By definition $\tau_6 - \tau_\gamma = \min \left( \tau_1 - \tau_2 - \epsilon_Y, \frac{1}{2} - \epsilon_W - \tau_\gamma \right)$. Taking $\epsilon_Y, \epsilon_W$ arbitrarily close to $0$, combining~\eqref{eqn: avgConv12} with~\eqref{eqn: avgConv1} and using the fact that $\tau_\gamma < \tau_1 - \tau_2$, we have that~\eqref{eqn: avgConv12} holds for every $0 \leq \tau_0 < \min\left( \tau_\gamma, \frac{1}{2} - \tau_\gamma \right)$. Moreover~\eqref{eqn: avgConv12} holds for every $\omega$ that satisfies Lemma~\ref{lem: staysBounded}, which yields~\eqref{eqn: averageConv}. 
\end{IEEEproof}

\subsection{Proof of Theorem~\ref{thm: main}}
\begin{IEEEproof}
	We show that, for \textit{all} $n \in \left\{1, \dots, N \right\}$, as long as the resilience condition~\eqref{eqn: resilienceCondition} (i.e., $\lambda_\min\left( \mathcal{G}_{\mathcal{N}} \right) > \Delta_{\mathcal{A}}$) is satisfied, we have \[\mathbb{P} \left( \lim_{t \rightarrow \infty} \left(t+1 \right)^{\tau_0} \left\lVert x_n(t) - \theta^* \right\rVert_2 = 0 \right) = 1\] for every $0 \leq \tau_0 < \min\left(\tau_\gamma, \frac{1}{2} - \tau_\gamma\right).$ By the triangle inequality, we have
\begin{equation}\label{eqn: main1}
	\left\lVert x_n(t) - \theta^* \right \rVert_2 \leq \left\lVert \widehat{\mathbf{x}}_t \right \rVert_2 + \left\lVert \overline{\mathbf{e}}_t \right \rVert_2.
\end{equation}
Lemma~\ref{lem: resilientConsensus} states that
\begin{equation}\label{eqn: main2}
	\mathbb{P} \left( \lim_{t \rightarrow \infty} (t+1)^{\tau_3} \left\lVert \widehat{\mathbf{x}}_t \right \rVert = 0 \right) = 1,
\end{equation}
for every $0 \leq \tau_3 < \tau_\gamma + \tau_1 - \tau_2$. Lemma~\ref{lem: averageConv} states that
\begin{equation}\label{eqn: main3}
	\mathbb{P} \left(\lim_{t \rightarrow \infty} (t+1)^{\tau_0} \left\lVert \overline{\mathbf{e}}_t \right \rVert_2 = 0 \right) = 1,
\end{equation}
for every $0 \leq \tau_0 < \min \left( \tau_\gamma, \frac{1}{2} - \tau_\gamma \right)$. Combining~\eqref{eqn: main1},~\eqref{eqn: main2}, and~\eqref{eqn: main3} yields the desired result. 
\end{IEEEproof}

\section{Resilience Analysis: Proof of Theorem~\ref{thm: resilience}}~\label{sect: resilience}
In this section, we prove Theorem~\ref{thm: resilience}, which relates the resilience condition $\lambda_{\min} \left( \mathcal{G}_{\mathcal{N}} \right) > \left\lvert \mathcal{A} \right \rvert$ to sparse observability. Theorem~\ref{thm: resilience} states that a necessary condition for \textbf{SAGE} to tolerate attacks on \textit{any} $s$ measurement streams is global $2s$-sparse observability. In addition, if the set of all unique measurement vectors $h_p$ is an orthogonal set, then global $2s$-sparse observability is sufficient for \textbf{SAGE} to tolerate attacks on any $s$ measurement streams.
\begin{IEEEproof}
	First, we show that if $\lambda_{\min} \left( \mathcal{G}_{\mathcal{N}} \right) > \left\lvert \mathcal{A} \right \rvert$ for all $\left\lvert \mathcal{A} \right \rvert \leq s$, then, the set of all measurement streams $\mathcal{P}$ is $2s$-sparse observable. We prove the contrapositive. Suppose that $\mathcal{P}$ is not $2s$-sparse observable. Then, there exists a set $\mathcal{X} \subseteq \mathcal{P}$ with cardinality $\left\lvert \mathcal{X} \right \rvert = 2s$ such that $\mathcal{G}_{\mathcal{P} \setminus \mathcal{X}} = \mathcal{H}^\intercal_{\mathcal{P}\setminus \mathcal{X}} \mathcal{H}_{\mathcal{P} \setminus \mathcal{X}}$ is not invertible. Partition $\mathcal{X}$ as $\mathcal{X} = \mathcal{A} \cup \mathcal{A}'$, where $\left \lvert \mathcal{A} \right \rvert = s$, $\left\lvert \mathcal{A'} \right \rvert = s$, and $\mathcal{A} \cap \mathcal{A}' = \emptyset$, and define $\mathcal{N} = \mathcal{P}\setminus \mathcal{A}$. Note that
\begin{align}\label{eqn: resProof1}
	\mathcal{G}_{\mathcal{N}} &= \sum_{p \in \mathcal{P}\setminus \mathcal{A}} h_p h_p^\intercal =\sum_{\mathcal{P}\setminus\mathcal{X}}  h_p h_p^\intercal + \sum_{\mathcal{A}'}  h_p h_p^\intercal \\
	&=\mathcal{G}_{\mathcal{P} \setminus \mathcal{X}}  + \mathcal{G}_{\mathcal{A}'},
\end{align}
The minimum eigenvalue of $\mathcal{G}_{\mathcal{N}}$ satisfies~\cite{Matrices2} \begin{equation}\label{eqn: rayleigh}\lambda_{\min}\left(\mathcal{G}_{\mathcal{N}}\right) \leq \min_{v \in \mathbb{R}^M, \: \left\lVert v \right\rVert_2 = 1} v^\intercal \mathcal{G}_{\mathcal{N}} v. \end{equation}
For any nonzero, unit norm $v$ in the nullspace of $\mathcal{G}_{\mathcal{P}\setminus \mathcal{X}}$, we have
\begin{equation}\label{eqn: resProof2}
	v^\intercal \mathcal{G}_{\mathcal{N}} v \leq \sum_{p \in \mathcal{A}'} v^\intercal h_p h_p^ \intercal v \leq \left\lvert \mathcal{A}' \right\rvert,
\end{equation}
where the last inequality in~\eqref{eqn: resProof2} follows as a consequence of the Cauchy-Schwarz Inequality~\cite{Matrices2} $h_p^\intercal v \leq \left\lVert h_p \right\rVert_2\left\lVert v \right\rVert_2 = 1$. From~\eqref{eqn: rayleigh} and~\eqref{eqn: resProof2}, we have that $\lambda_{\min} \left(\mathcal{G}_{\mathcal{N}} \right) \leq \left\lvert \mathcal{A} \right\rvert$, which shows the contrapositive. 

Second, we show that, under the condition that $\mathcal{V}$, the set of \textit{unique} measurement vectors $h_p$, is orthogonal, if $\mathcal{P}$ is $2s$-sparse observable, then, $\lambda_{\min} \left( \mathcal{G}_{\mathcal{N}} \right) > \left\lvert \mathcal{A} \right \rvert$ for all $\left\lvert \mathcal{A} \right \rvert \leq s$. 
We resort to contradiction. Let $\mathcal{V} = \left\{v_1, \dots, v_{\left\lvert \mathcal{V} \right\rvert} \right\}$. If $\left\lvert \mathcal{V} \right\rvert < M$, then, define $v_{\left\lvert \mathcal{V} \right \rvert + 1}, \dots, v_M$ so that $\left\{ v_1, \dots, v_M \right\}$ forms an orthonormal basis for $\mathbb{R}^M$. Suppose there exists $\mathcal{A}$ such that $\left\lvert \mathcal{A} \right\rvert \leq s$ and $\lambda_{\min} \left(\mathcal{G}_{\mathcal{N}} \right) \leq \left\lvert \mathcal{A} \right \rvert$. Define, for $m = 1, \dots M$, 
\begin{equation}\label{eqn: resProof4}
	\mathcal{N}_m = \left\{ p \in \mathcal{N} | h_p = v_m \right\}.
\end{equation}
Then, for any $v_m$, $m = 1, \dots, M$, we have
\begin{equation}\label{eqn: resProof4a}
	\mathcal{G}_{\mathcal{N}} v_m = \sum_{p \in \mathcal{N}} h_p h_p^\intercal v_m = \left\lvert \mathcal{N}_m \right\rvert v_m,
\end{equation}
since, by the orthogonality of $\mathcal{V}$, we have $h_p^\intercal v_m = 1$ if $h_p = v_m$ and $h_p^\intercal v_m = 0$ if $h_p \neq v_m$. Equation~\eqref{eqn: resProof4a} states that the eigenvalues of $\mathcal{G}_{\mathcal{N}}$ are $\left\lvert \mathcal{N}_1 \right \rvert, \dots, \left\lvert \mathcal{N}_M \right \rvert$. Let $m^* = \argmin_{m \in \{1, \dots, M\}} \left\lvert \mathcal{N}_m \right \rvert,$ and define $\mathcal{X} = \mathcal{A} \cup \mathcal{N}_{m^*}.$ Since $\lambda_{\min} \left(\mathcal{G}_{\mathcal{N}} \right) \leq \left\lvert \mathcal{A} \right\rvert \leq s$, we have $\left\lvert \mathcal{X} \right \rvert \leq 2s$. Then, we calculate the minimum eigenvalue of $\mathcal{G}_{\mathcal{P} \setminus \mathcal{X}}$ as
\begin{equation}\label{eqn: resProof5}
	\lambda_{\min} \left( \mathcal{G}_{\mathcal{P} \setminus \mathcal{X}} \right) = \lambda_{\min} \left( \mathcal{G}_{\mathcal{N} \setminus \mathcal{N}_{m^*}} \right) = 0.
\end{equation}
Equation~\eqref{eqn: resProof5} indicates that the set of all measurement streams $\mathcal{P}$ is not $2s$-sparse observable, which is a contradiction and shows the desired result. 
\end{IEEEproof}


\section{Numerical Examples}\label{sect: examples}
We demonstrate the performance of \textbf{SAGE} through two numerical examples. 

\subsection{Homogeneous Measurements}
In this first example, we consider a random geometric network of $N = 500$ agents each making (homogeneous) measurements of an unknown parameter $\theta^* \in \mathbb{R}^2$. 
\begin{figure}[h!]
	\centering
	\includegraphics[width = 0.7\columnwidth]{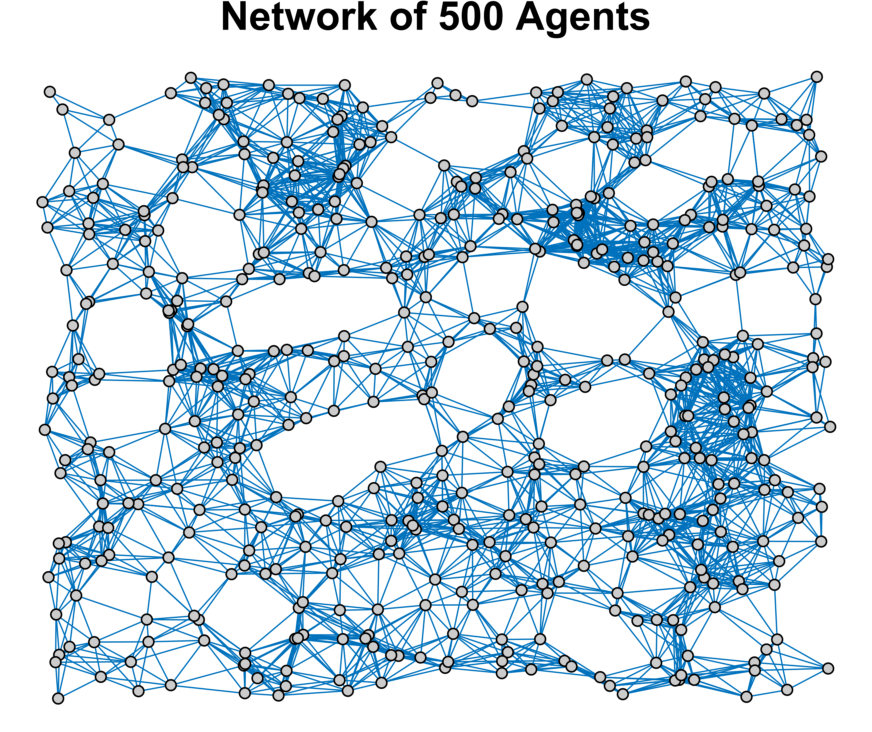}
	\caption{A random geometric network of $N = 500$ agents.}\label{fig: network1}
\end{figure}
The parameter $\theta^*$ may represent, for example, the location of a target to be tracked. Each agent makes a local stream of measurements $y_n(t) = \theta^* + w_n(t)$, where the measurement noise $w_n(t)$ is additive white Gaussian noise (AWGN) with mean $0$ and covariance $\sigma^2 I_2$. {\color{black} According to~\eqref{eqn: relaxedResilience}, \textbf{SAGE} is resilient to attacks on \textit{any} set of $250$ measurement streams. }

In this example, we consider four levels of local signal-to-noise ratio (SNR): -7 dB, -13 dB, -19 dB, and -25 dB. Note that, at these SNRs, even in the absence of attack, the noise at each local device is much stronger than the measured signal. We place agents uniformly at random on a two dimensional grid, and agents may communicate with nearby agents whose Euclidean distance is below a certain threshold. Figure~\ref{fig: network1} depicts the network setup for the first numerical example. {\color{black} Communication links may fail at random. In each iteration, each communication link fails (independently) with probability $0.1$ (i.e., $10\%$ chance of failure in each iteration). } An adversary compromises the measurement streams of $100$ of the agents, chosen uniformly at random. In this example, the agents under attack have measurement streams $y_n^a(t) = -3 \theta^* + w_n(t)$. {\color{black} \footnote{\color{black} Although we use a fixed attack $a_n(t)$ in this simulation, recall that, in general, the attack $a_n(t)$ may be unbounded and time varying.}}

{\color{black}It is important to correctly choose $\alpha_t = \frac{a}{(t+1)^{\tau_1}}$, $\beta_t = \frac{b}{(t+1)^{\tau_2}}$, and $\gamma_t = \frac{\Gamma}{(t+1)^{\tau_\gamma}}$. The weight selection procedure (equations~\eqref{eqn: alphaBeta} and~\eqref{eqn: gammaDef}) only require $a, b, \Gamma > 0$, $0 < \tau_2 < \tau_1 \leq 1$, and $\tau_\gamma <\min \left( \frac{1}{2}, \tau_1 - \tau_2 \right)$. We use the following guidelines for selecting these weights in our numerical simulations. We choose
\[a = 1, b = \frac{1}{\lambda_N \left(L \right)},\] 
where $L$ is the Laplacian of the fixed graph shown in Figure~\ref{fig: network1} (with no failing links). The effect of choosing different $\Gamma$ on the performance of \textbf{SAGE} also depends on the true value of $\theta^*$ and the adversary's actions $a_n(t)$. Since we do not know these values, instead of providing a rule for choosing $\Gamma$, we show, for fixed $\theta^*$ and $a_n(t)$, the effect of choosing different $\Gamma$ on the estimator's performance. For the decay rates $\tau_1, \tau_2$, and $\tau_\gamma$, we choose
\[\tau_1 = 0.26, \tau_2 = 0.001, \tau_\gamma = 0.25. \]

We now explain our guidelines. From~\eqref{eqn: staysBounded7}, we see that the evolution of $\left\lVert \overline{\mathbf{e}}_t \right\rVert_2$ depends on $\left \lVert I_{M} - \frac{\alpha_t}{N} \sum_{p \in \mathcal{N}} k_p(t) h_p h_p^\intercal \right \rVert_2$. If we choose $a = 1$, then, we guarantee that \[\left \lVert I_{M} - \frac{\alpha_t}{N} \sum_{p \in \mathcal{N}} k_p(t) h_p h_p^\intercal \right \rVert_2 \leq 1,\] for all $t$, which ensures, at each time step, a contractive effect on the evolution of $\left\lVert \overline{\mathbf{e}}_t \right\rVert_2$. Similarly, from~\eqref{eqn: resCon2}, we see that the evolution of $\left\lVert \widehat{\mathbf{x}}_t \right\rVert_2$ depends on 
$\left\lVert I_{NM}\! -\! \beta_t L_t \!\otimes\!I_M \!-\! P_{NM}\right\rVert_2.$ If we consider the fixed graph Laplacian $L$ and choose $b = \frac{1}{ \lambda_{N} \left( L \right)}$, then, we guarantee that
\[ \left\lVert  I_{NM}\! -\! \beta_t L_t \!\otimes\!I_M \!-\! P_{NM} \right\rVert_2 \leq 1,\]
for all $t$, which ensures, at each time step, a contractive effect on the evolution of $\left\lVert \widehat{\mathbf{x}} \right\rVert_2$. For choosing the decay rates, note that, in Theorem~\ref{thm: main}, the local estimates' convergence rate is upper bounded by $\min\left(\tau_\gamma, \frac{1}{2}-\tau_\gamma \right)$. We maximize this upper bound by selecting $\tau_\gamma = 0.25$. To promote fast resilient consensus, we select a small decay rate for the consensus weight ($\tau_2 = 0.001$). Finally, we choose $\tau_1$ to be large enough ($\tau_1 = 0.26$) to satisfy $\tau_\gamma < \tau_1 - \tau_2$. Choosing $\tau_1$ to be too large causes the innovation weight $\alpha_t$ to decay rapidly, which may slow down the convergence of the algorithm.}

The agents follow \textbf{SAGE} with the following weights: $a = 1, \tau_1 = 0.26, b = 0.0337, \tau_2 = 0.001, \Gamma = 5, \tau_\gamma = 0.25$. We compare the performance of \textbf{SAGE} against the performance of the baseline consensus+innovations estimator~\cite{Kar1, Kar2} with the same weights $a, \tau_1, b, \tau_2$. 

\begin{figure}[h!]
	\centering
	\includegraphics[width = 0.85\columnwidth]{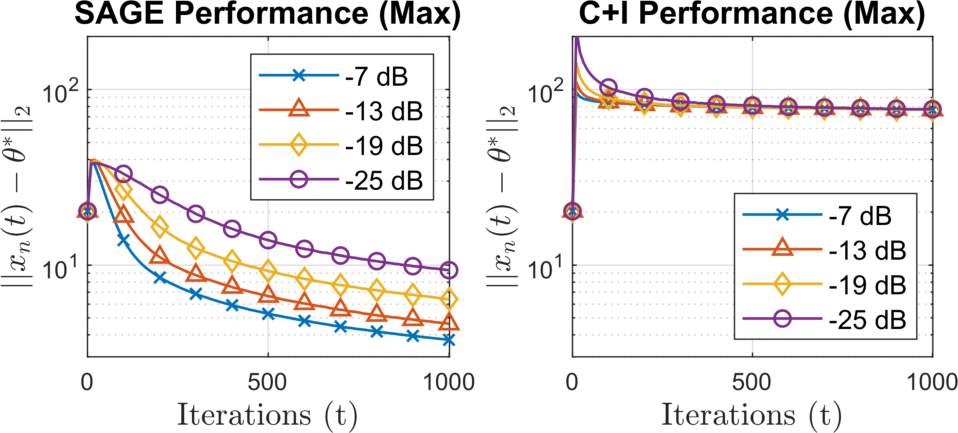}
	\caption{Maximum over all agents of the local estimation root mean square error (RMSE) for \textbf{SAGE} (left) and baseline consensus+innovations (right).}\label{fig: SAGEMax}
\end{figure}


Figures~\ref{fig: SAGEMax} shows the evolution of the maximum of the local root mean squared error (RMSE) (across all agents) over iterations of \textbf{SAGE} and the baseline consensus+innovations estimator. We report the average of the maximum RMSE over $500$ trials. In the presence of adversaries, under \textbf{SAGE}, the local RMSE decreases with increasing number of iterations, while, under baseline consensus+innovations, there is a persistent local RMSE that does not decrease. Figure~\ref{fig: SAGEMax} shows that, under \textbf{SAGE}, higher local SNR yields lower RMSE. 

{\color{black}
\begin{figure}[h!]
	\centering
	\includegraphics[width = 0.85\columnwidth]{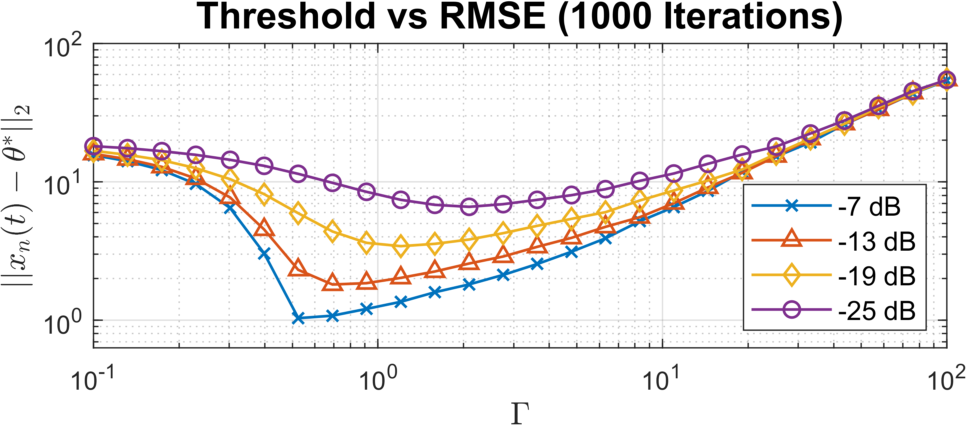}
	\caption{Maximum local RMSE after $1000$ iterations versus $\Gamma$.}\label{fig: SAGEGamma}
\end{figure}
Figure~\ref{fig: SAGEGamma} shows the effect of choosing different $\Gamma$ on the maximum local RMSE (across all aents) after $1000$ iterations of \textbf{SAGE}. If we choose $\Gamma$ to be too small (e.g., at $-7 \text{ dB}$ SNR, if $\Gamma < 10^{-0.8}$), then, the thresholding is too aggressive, and we limit the contributions of noncompromised measurements to the estimate update. As we increase $\Gamma$, \textbf{SAGE} incorporates more information from the noncompromised measurements, resulting in lower RMSE. If we increase $\Gamma$ too much, however, \textbf{SAGE} becomes less effective at limiting the impact of the attacks, which results in higher RMSE. The ``optimal'' choice of $\Gamma$ (the choice of $\Gamma$ that provides the lowest RMSE in Figure~\ref{fig: SAGEGamma}) depends on the value of the parameter $\theta^*$ and the attack $a_n(t)$, which we do not know a priori. }

The performance of \textbf{SAGE} also depends on the number of compromised measurement streams. 
\begin{figure}[h!]
	\centering
	\includegraphics[width = 0.85\columnwidth]{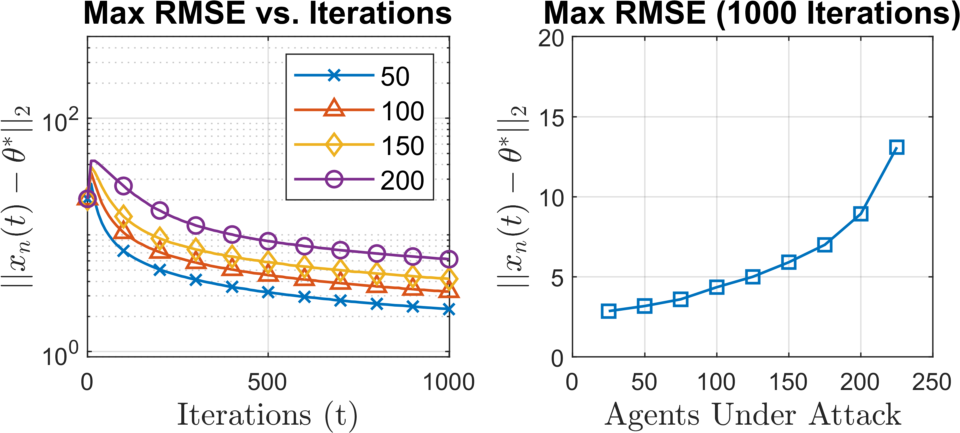}
	\caption{Maximum local estimation RMSE versus number of \textbf{SAGE} iterations (left). Maximum local RMSE after $1000$ iterations versus number of agents with compromised measurements (right).}\label{fig: SAGENumAttacked}
\end{figure}
Figure~\ref{fig: SAGENumAttacked} describes the relationship between the number of agents with compromised measurements and the maximum local RMSE (at a local SNR of $-13$ dB). We report the average maximum RMSE over $500$ trials, where in each trial, the adversary attacks a new set of agents chosen uniformly at random. Figure~\ref{fig: SAGENumAttacked} shows that, for a fixed number of iterations of \textbf{SAGE}, the local RMSE increases as the number of agents under attack increases.{\color{black}\footnote{\color{black} When there is no measurement noise and the parameter is guarnteed to be bounded (i.e., $\left\lVert \theta^* \right\rVert \leq \eta$ for some $\eta > 0$), reference~\cite{ChenDistributed2} provides an analytical description of the relationship between the number of agents under attack and the local RMSE.}} {\color{black} We fit the exponential curve
\begin{equation*}
	\text{RMSE}_{1000} = 1.756 \times e^{0.00855 \left\lvert \mathcal{A} \right\rvert},
\end{equation*}
with $R^2 = 0.9603$, to describe the relationship between the number of compromised measurements $\left\lvert \mathcal{A} \right\rvert$ and the RMSE after $1000$ iterations of $\textbf{SAGE}$.} 

\subsection{Heterogeneous Measurements}
\begin{figure}[h!]
	\centering
	\includegraphics[width = 0.7\columnwidth]{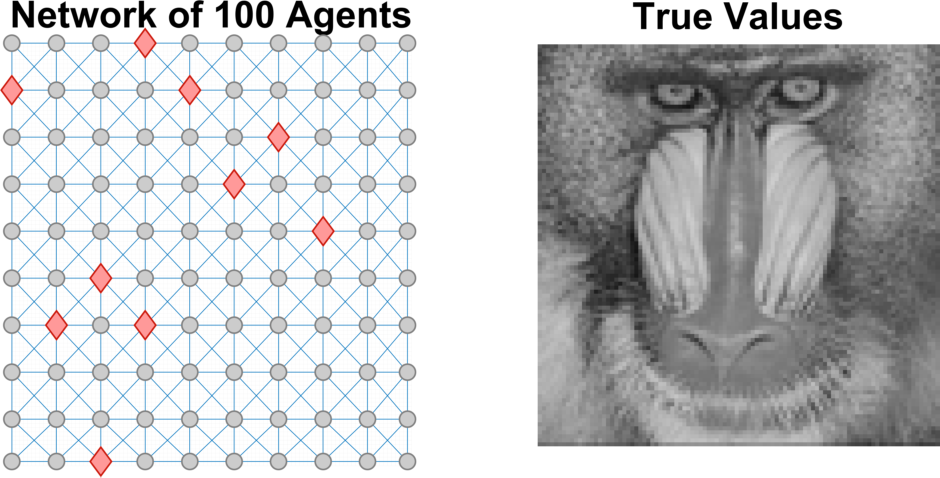}
	\caption{\color{black}(Left) A network of $N = 100$ robots. An adversary compromises the measurement streams of $10$ robots (shown as red diamonds). (Right) The true value of $\theta^*$, where each pixel value of the image represents the safety score at that particular grid location.}\label{fig: network2}
\end{figure}

In the second example, we demonstrate the performance of \textbf{SAGE} with heterogeneous measurement models. The physical motivation for the example is as follows. {\color{black} We consider the network of $N = 100$ robots placed in a two dimensional environment modeled by a $100 \text{ unit } \times 100 \text{ unit }$ grid (see Figure~\ref{fig: network2}). Similar to the first example, the communication links in the network fail at random. In each iteration, each link fails independently with probability $0.1$.} Each square of the grid has an associated $[0, 255]$-valued safety score, which represents the state of that particular location. For example, a safety score of $0$ may mean that there is an inpassable obstacle at a particular location, while a safety score of safety score of $255$ may represent a location that is completely free of obstacles. {\color{black} The unknown parameter $\theta^* \in [0, 255]^{10,000}$ is the collection of safety scores over the entire environment. Each component of $\theta^*$ is the safety score at a single grid location.} We represent the true value of $\theta^*$ an image (the image of the baboon in Figure~\ref{fig: network2}), where each pixel value represents the score of a grid location. 

The agents' goal is to recover the \textit{entire} image from their local measurements. Each agent measures the local pixel values in a 45 pixel by 45 pixel grid centered at its own location. {\color{black} The measurement of agent $n$ is 
\begin{equation*}
	y_n(t) = \underbrace{\left[\begin{array}{c} h_{\overline{P}_n + 1}^\intercal \\ \vdots \\ h_{\overline{P}_n + P_n} \end{array} \right]}_{H_n} \theta^* + w_n(t),
\end{equation*}
where each measurement vector $h_p$ is a canonical basis vector. In this example, the $H_n$ matrices are chosen to emphasize local measurements. The measurement matrix $H_n$ ``selects'' in this simulation the components of $\theta^*$ corresponding to all grid locations within a $45$ pixel by $45$ pixel subgrid centered at the location of agent $n$. The measurement noise $w_n(t)$ is i.i.d. $\mathcal{N} \left(0,  100 I_{P_n}\right).$ The dimension $P_n$ of the measurement $y_n(t)$ depends on the location of the agent $n$. Agents in the center of the grid environment measure the states of $P_n = 45 \times 45 = 2025$ grid locations, while agents toward the edge of the environment have smaller $P_n$ (for example, agents in the corner of the environment have $P_n = 22 \times 22 = 484$). An adversary attacks $10$ agents and compromises all of their measurement streams: all compromised measurement streams $p \in \mathcal{A}$ take value $y_n^{(p)}(t) = 255$. In total, there are $P = 93,025$ measurement streams, and the adversary  attacks $10$ agents and compromises $\left\lvert \mathcal{A} \right\rvert = 16,082$ measurement streams.\footnote{\color{black} The set $\mathcal{A}$ is the set of compromised measurement streams, not the set of compromised agents.} It is impractical to explicitly detect and identify each compromised measurement stream in this scenario, due to the combinatorial cost.

According to~\eqref{eqn: relaxedResilience} and Theorem~\ref{thm: resilience}, \textbf{SAGE} is guaranteed to be resilient to attacks on \textit{any} $4$ measurement streams, which is the same as the most resilient \textit{centralized} estimator. That is, there exists a specific set of $5$ measurement streams, which, if compromised, prevents \textbf{SAGE} (and any other estimator) from consistently estimating $\theta^*$. Even though there exists a \textit{specific} set of $5$ (or more) compromised measurement streams prevent \textbf{SAGE} from producing consistent estimates, in this simulation, we show that \textbf{SAGE} may still produce consistent estimates even the number of compromised measurements vastly exceeds four   ($\left\lvert \mathcal{A} \right\rvert = 16,082 \gg 4$).  }

{\color{black}We use the same guidelines as in the first example to select the weights for \textbf{SAGE}:} $a = 1, \tau_1 = 0.26, b = 0.0494, \tau_2 = 0.001, \Gamma = 100, \tau_\gamma = 0.25$. We compare the performance of \textbf{SAGE} against the performance of the baseline consensus+innovations estimator.
\begin{figure}[h!]
	\centering
	\includegraphics[width = 0.85\columnwidth]{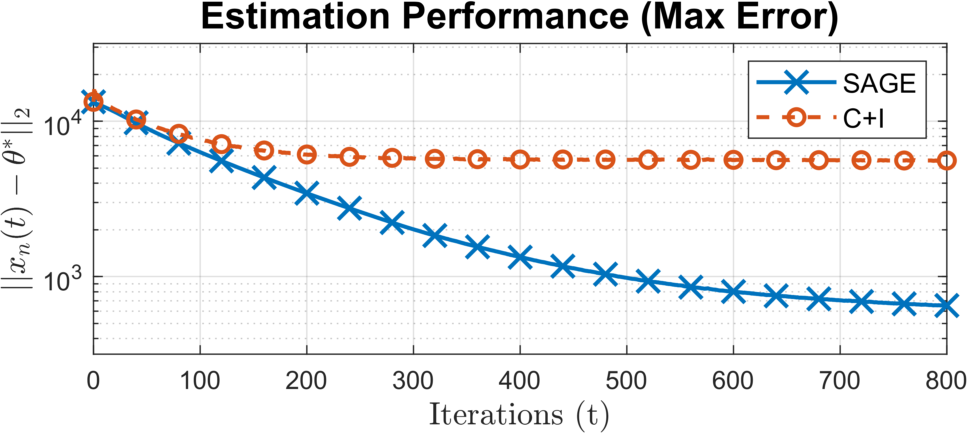}
	\caption{\color{black} Maximum over all agents of the local estimation root mean square error (RMSE) for \textbf{SAGE} and baseline consensus+innovations.}\label{fig: SAGEImageIteration}
\end{figure}
\begin{figure}[h!]
	\centering
	\includegraphics[width = 0.7\columnwidth]{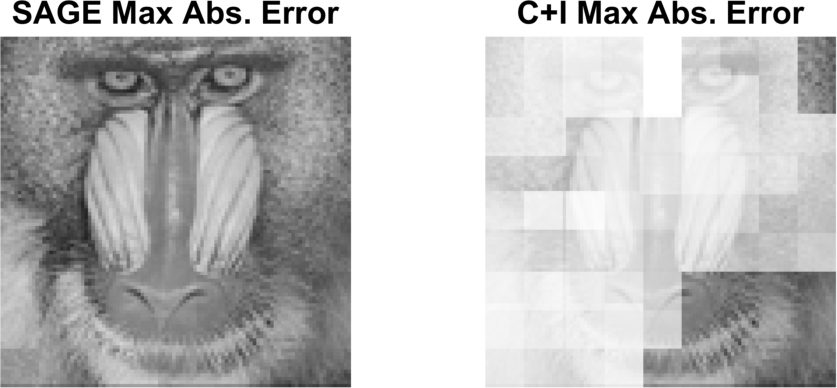}
	\caption{\color{black} Comparison the estimation result after $800$ iterations of \textbf{SAGE} (left) and baseline consensus+innovations (right).}\label{fig: SAGEComparison}
\end{figure}

{\color{black} Figure~\ref{fig: SAGEImageIteration} shows the evolution of the maximum local RMSE (over all agents) for \textbf{SAGE} and the baseline consensus+innovations estimator. We report the average of the maximum RMSE over $500$ trials. Under the measurement attacks, \textbf{SAGE} produces estimates with decreasing RMSE (with increasing number of iterations), while the baseline consensus+innovations estimator has a persistent local RMSE that does not decrease.} Figure~\ref{fig: SAGEComparison} compares the estimation results of \textbf{SAGE} and the baseline consensus+innovations estimator, after $800$ iterations of each algorithm. We reconstruct each pixel of the image using the worst estimate across all agents. Figure~\ref{fig: SAGEComparison} shows that, under \textbf{SAGE}, \textit{all} of the agents resiliently recover the image even when there are measurement attacks. In contrast, under the baseline consensus+innovations estimator, measurement attacks prevent the agents from consistently estimating $\theta^*$. 

\section{Conclusion}\label{sect: conclusion}
In this paper, we have studied resilient distributed estimation under measurement attacks. A network of IoT devices makes heterogeneous, linear, successive noisy measurements of an unknown parameter $\theta^*$ and cooperates over  a sparse communication network to sequentially process their measurements and estimate the value of $\theta^*$. An adversary attacks a subset of the measurements and manipulates their values arbitrarily. This paper presented \textbf{SAGE}, the Saturating Adaptive Gain Estimator, a recursive, consensus+innovations estimator that is resilient to measurement attacks. 

Under \textbf{SAGE}, each device applies an adaptive gain to its innovation (the difference between its observed and predicted measurements) to ensure that the magnitude of the scaled each innovation component is below time-varying, decaying threshold. This adaptive gain limits the impact of compromised measurements. As long as the number of compromised measurements is below a particular bound, then, we demonstrate that, for any on-average connected topology, \textbf{SAGE} guarantees the strong consistency of all of the (compromised and uncompromised) devices' local estimates.  When each measurement stream collects data about a single component of $\theta^*$, \textbf{SAGE} achieves the same level of resilience (in terms of number of tolerable compromised measurement streams) as the most resilient centralized estimator. For example, this could occur in air quality monitoring, when devices measure local pollutant concentrations corresponding to inidivudal components of the unknown parameter. Finally, we illustrated the performance of \textbf{SAGE} through numerical examples. {\color{black} Future work includes designing resilient estimators for dynamic (time-varying) parameters and dealing with measurement attacks that roam over the network.}

\appendix
The proof of Theorem~\ref{thm: main} requires several intermediate results. The following result from~\cite{ChenDistributed1} characterizes the behavior of time-averaged measurement noise.
\begin{lemma}[Lemma 5 in~\cite{ChenDistributed1}]\label{lem: measureNoise}
	Let $v_1, v_2, v_3, \dots$ be i.i.d. random variables with mean $\mathbb{E} \left[ v_t \right] = 0$ and finite covariance $\mathbb{E} \left[ v_t v_t^{\intercal} \right] = \Sigma.$ Define the time-averaged mean
	$\overline{v}_t = \frac{1}{t+1} \sum_{j=0}^t v_t.$
Then, we have
\begin{equation}\label{eqn: noiseConvergence}
	\mathbb{P} \left( \lim_{t \rightarrow \infty} (t+1)^{\delta_0} \left\lVert \overline{v}_t \right\rVert_2 = 0 \right) = 1,
\end{equation}
for all $0 \leq \delta_0 < \frac{1}{2}.$
\end{lemma}

We will need to study the convergence properties of scalar, time-varying dynamical systems of the form:
\begin{equation}\label{eqn: timeVaryingSystem}
	w_{t+1} = \left( 1 - r_1(t) \right) w_t + r_2(t),
\end{equation}
where
\begin{equation}\label{eqn: decayRateConditions}
	r_1(t) = \frac{c_1}{(t+1)^{\delta_1}}, r_2(t) = \frac{c_2}{(t+1)^{\delta_2}},
\end{equation}
$c_1, c_2 > 0$, and $0 < \delta_1 < \delta_2 < 1$. The following Lemma from~\cite{Kar1} describes the convergence rate of the system in~\eqref{eqn: timeVaryingSystem}. 
\begin{lemma}[Lemma 5 in~\cite{Kar1}]\label{lem: timeVaryingConvergence}
	The system in~\eqref{eqn: timeVaryingSystem} satisfies
\begin{equation}
	\lim_{t \rightarrow \infty} \left(t+1\right)^{\delta_0} w_{t+1} = 0,
\end{equation}
for every $0 \leq \delta_0 < \delta_2 - \delta_1.$
\end{lemma}


We will also need the convergence properties of the system in~\eqref{eqn: timeVaryingSystem} when the $r_1(t)$ is random.

\begin{lemma} [Lemma 4.2 in~\cite{Kar4}]\label{lem: timeVaryingStochasticConvergence}
	Let $\left\{w_t \right\}$ be an $\mathcal{F}_t$-adapted process satisfying~\eqref{eqn: timeVaryingSystem}, i.e., $ w_{t+1} = \left( 1 - r_1(t) \right) w_t + r_2(t),$ where $\left\{r_1(t) \right\}$ is an $\mathcal{F}_{t+1}$-adapted process such that, for all $t$, $0 \leq r_1(t) \leq 1$ a.s. and
\begin{equation}\label{eqn: randomConvergenceCondition}
	\frac{c_1}{(t+1)^{\delta_1}} \leq \mathbb{E} \left[ r_1(t) \vert \mathcal{F}_t \right] \leq 1,
\end{equation}
with $c_1 > 0$, $0 \leq \delta_1 < 1$.  Let $r_2 = \frac{c_2}{(t+1)^{\delta_2}}$ with $c_2 > 0$ and $\delta_1 < \delta_2 < 1$. Then, we have
\begin{equation}
	\mathbb{P} \left( \lim_{t \rightarrow \infty} (t+1)^{\delta_0} w_t = 0 \right) = 1,
\end{equation}
for every $0 \leq \delta_0 < \delta_2 - \delta_1$. 
\end{lemma}

To account for the effect of random, time-varying Laplacians, we use the following result from~\cite{Kar4}.
\begin{lemma} [Lemma 4.4 in~\cite{Kar4}]\label{lem: randomGraphWeights}
	Let $\mathcal{C} \subseteq \mathbb{R}^{NM}$ be the consensus subspace, defined as
\begin{equation}
	\mathcal{C} = \left\{ w \in \mathbb{R}^{NM} \vert w = \mathbf{1}_N \otimes v, v \in \mathbb{R}^N \right\},
\end{equation}
and let $\mathcal{C}^\perp$ be the orthogonal complement of $\mathcal{C}$. Let $\left\{ w_t \right\}$ be an $\mathcal{F}_t$-adapted process such that $w_t \in \mathcal{C}^{\perp}$. Also, let $\left\{L_t \right\}$ be an i.i.d. sequence of Laplacian matrices that is $\mathcal{F}_{t+1}$-adapted, independent of $\mathcal{F}_t$, and satisfies $\lambda_2 \left( \mathbb{E} \left[ L_t \right] \right) > 0.$ Then, there exists a measurable $\mathcal{F}_{t+1}$-adapted process $\left\{r(t) \right\}$ and a constant $c_r > 0$ such that $0 \leq r(t) \leq 1$ a.s. and 
\begin{equation}
 	 \left\lVert  \left(I_{NM} - \beta_t L_t \otimes I_M \right) w_t \right\rVert_2 \leq \left(1- r(t) \right) \left\lVert w_t \right \rVert_2,
\end{equation}
with
\begin{equation}
	\mathbb{E} \left[ r(t) \vert \mathcal{F}_t \right] \geq \frac{c_r}{(t+1)^{\tau_2}},
\end{equation}
for all sufficiently large $t$. The terms $\beta_t$ and $\tau_2$ are defined according to~\eqref{eqn: alphaBeta}.
\end{lemma}

We will also need to analyze scalar, time-varying dynamical systems of the form:
\begin{equation}\label{eqn: supSystem}
\begin{split}
	\widehat{m}_{t+1} &= \left(1 - \frac{r_1(t)}{\left( m_t + c_3 \right)} \right) m_t + r_2(t), \\
	m_{t+1} &= \max \left( \left \lvert \widehat{m}_{t+1} \right\rvert, \left\lvert m_t \right\rvert \right),
\end{split}
\end{equation}
with initial condition $m_0 \geq 0$, where $r_1(t)$ and $r_2(t)$ are given by~\eqref{eqn: decayRateConditions}, and $c_3 > 0$. 
\begin{lemma}\label{lem: supBounded}
	The system in~\eqref{eqn: supSystem} satisfies
	\begin{equation}
		\sup_{t \geq 0} m_t < \infty.
	\end{equation}
\end{lemma}
\begin{IEEEproof}[Proof of Lemma~\ref{lem: supBounded}]
	First, we show that there exists finite $T$ such that $m_{T+1} = m_{T}$ and $r_1(T) < 1$. We separately consider the cases of $T = 0$ and $T > 0$. If $m_1 = m_0$ and $r_1(0) < 1$, then the condition is satisfied for $T = 0$. Otherwise, by construction, $m_t$ is a non-negative, non-decreasing sequence.  Since $r_1(t)$ decreases in t, for $T$ large enough, we have $ r_1(T) < 1$ and $\left(1 - \frac{r_1(T)}{m_T + c_3} \right) \geq 0$, which means that $\left \lvert \widehat{m}_{T+1} \right \rvert = \widehat{m}_{T+1}$. Further, this means that, for $T$ large enough, a sufficient condition for $m_{T+1} = m_{T}$ is $m_{T} \geq \widehat{m}_{T+1}$. Computing $m_T - \widehat{m}_{T+1}$, we have
\begin{equation}\label{eqn: supProof1}
	m_T - \widehat{m}_{T+1} = \frac{c_1 m_T}{(T+1)^{\delta_1} \left(m_T + c_3\right)} - \frac{c_2}{(T+1)^{\delta_2}}.
\end{equation}
The term $\frac{m_T}{m_T + c_3}$ is increasing in $m_T$. Since $m_T$ is a nondecreasing sequence, we have
\begin{equation}\label{eqn: supProof2}
	\frac{m_T}{m_T+c_3} \geq \frac{m_1}{m_1 + c_3} > 0
\end{equation}
for all $T > 0$. To derive~\eqref{eqn: supProof2} from~\eqref{eqn: supProof1}, we have used the fact that $m_1 > 0$, which is guaranteed by~\eqref{eqn: supSystem}. 

Then, letting $c_4 = \frac{m_1}{m_1 + c_3}$, substituting into~\eqref{eqn: supProof1}, and performing algebraic manipulations, a sufficient condition for $m_t \geq \widehat{m}_{T+1}$ is
\begin{equation}\label{eqn: supProof3}
	T \geq  \left( \frac{c_2}{c_1 c_4} \right)^{\frac{1}{\delta_2 - \delta_1}} - 1,
\end{equation}
which shows that there exists finite $T$ such that $m_{T+1} = m_{T}$. Second, we show that $m_t = m_T$ for all $t \geq T$. If $t \geq T$, then $t$ also satisifies the sufficient condition in~\eqref{eqn: supProof3}, which means that $m_T = m_{T+1} = m_{T+2} = \dots$. Finally, we have
$\sup_{t \geq 0} m_t = m_T < \infty.$
\end{IEEEproof}
\bibliography{IEEEabrv,References}

\end{document}